\documentclass[12pt]{amsart}
\usepackage{amsmath, amsthm, amssymb, amsfonts}
\usepackage{mathtools}
\usepackage{mathrsfs, changes} 
\usepackage{dsfont}
\usepackage{hyperref}\hypersetup{
    colorlinks=true,
    linkcolor=purple,
    filecolor=blue,      
    urlcolor=cyan,
    }
\usepackage{comment}
\usepackage{enumitem}

\numberwithin{equation}{section}  %numberwithin goes before cleverefs when using hyperref
\usepackage[sort&compress,capitalize,nameinlink]{cleveref}
\crefname{app}{Appendix}{Appendices}

\crefrangeformat{equation}{\upshape(#3#1#4)\textendash(#5#2#6)}

\usepackage{tikz, pgfplots}
\usepgfplotslibrary{statistics}
\usetikzlibrary{decorations.markings,calc,math}
\pgfplotsset{compat=1.18}
\usepackage{subcaption}
\usepackage{diagbox}
\usepackage{wrapfig}
\usepackage{overarrows}
\pgfmathdeclarefunction{gamma}{1}{\pgfmathparse{
    (2.506628274631*sqrt(1/#1) + 0.20888568*(1/#1)^(1.5) + 0.00870357*(1/#1)^(2.5) - (174.2106599*(1/#1)^(3.5))/25920 - (715.6423511*(1/#1)^(4.5))/1244160)*exp((-ln(1/#1)-1)*#1)}}
\pgfmathdeclarefunction{poisson_pmf}{2}{%
\pgfmathparse{(#2^#1)*exp(-#2)/(gamma(#1+1))}%
}
\pgfmathdeclarefunction{normal_pdf}{3}{%
\pgfmathparse{1/(#3*sqrt(2*pi))*exp(-((#1-#2)^2)/(2*#3^2))}%
}

\theoremstyle{plain}
\newtheorem{theorem}{Theorem}

\newtheorem*{corollary*}{Corollary}
\newtheorem{lemma}[theorem]{Lemma}

\theoremstyle{definition}
\newtheorem{definition}[theorem]{Definition}
\newtheorem*{definition*}{Definition}
\newtheorem{assumption}[theorem]{Assumption}

\newtheorem*{hypothesis*}{Hypothesis}
\theoremstyle{remark}
\newtheorem{remark}{Remark}
\newtheorem*{remark*}{Remark}
\newtheorem{example}{Example}
\numberwithin{theorem}{section}
\numberwithin{remark}{section}
\numberwithin{example}{section}

\newcommand{\Expect}{\mathbb{E}}
\newcommand{\Prob}{\mathsf{\mathbb{P}}}

\newcommand{\R}{\mathbb{R}}

\newcommand{\Z}{\mathbb{Z}}
\newcommand{\N}{\mathbb{N}}

\DeclarePairedDelimiter{\braces}{\{}{\}}
\DeclarePairedDelimiter{\parens}{(}{)}

\DeclarePairedDelimiter{\floor}{\lfloor}{\rfloor}
\newcommand{\zero}{\boldsymbol{0}}
\DeclareMathOperator{\pne}{\mathsf{PNE}}
\DeclareMathOperator{\Abso}{\mathsf{Abso}}
\DeclareMathOperator{\des}{\mathsf{des}}
\DeclareMathOperator{\Des}{\mathcal{D}}

\newcommand{\vertices}{\mathcal{V}_{N}}
\newcommand{\edges}{\mathcal{E}_{N}}
\NewOverArrowCommand{longvec}{}
\newcommand{\oedges}{\longvec*{\mathcal{E}}_{N}}
\newcommand{\orient}[1]{\overrightarrow{#1}}
\newcommand{\rorient}[1]{\overleftarrow{#1}}

\newcommand{\card}[1]{\mathsf{{card}}\parens*{#1}}
\newcommand{\ball}{{\rm Ball}}

\newcommand{\cube}{\mathcal{H}_{N}}
\newcommand{\ocube}{\orient{\mathcal{H}}_{N}^{\beta}}

\newcommand{\ularg}{\mathcal{L}_{N}^{\beta}}
\newcommand{\olargo}{\orient{\mathcal{L}}_{N}^{{\beta},{\zero}}}
\newcommand{\rlargo}{\rorient{\mathcal{L}}_{N}^{{\beta},{\zero}}}
\newcommand{\trap}{\mathscr{T}}

\title[Finding Pure Nash Equilibria in large Random Games]{Finding Pure Nash Equilibria in large Random Games}

\author{Andrea Collevecchio}
\address{School of Mathematics and Monash Data Futures Institute,
Monash  University, Melbourne, Australia}
\email{Andrea.Collevecchio@monash.edu}
\author{Tuan-Minh Nguyen}
\address{School of Mathematics, Monash  University, Melbourne, Australia}
\email{tuanminh.nguyen@monash.edu}
\author{Ziwen Zhong}
\address{School of Mathematics, Monash  University, Melbourne, Australia}
\email{ziwen.zhong@monash.edu}

\begin{document}
\begin{abstract} 
Best Response Dynamics (BRD) constitute a family of strategy update rules designed to identify Pure Nash Equilibria (PNEs) in strategic games. In each iteration, a player is selected uniformly at random and adopts a best response strategy, maximising their payoff given the current strategies of the other players. If no improvement is possible for the selected player, another is chosen. The process terminates upon reaching a state where no player wishes to deviate, corresponding to a PNE. However, convergence is not guaranteed {as either PNEs may not exist, or the dynamics may become confined within substructures of the game, which we call traps, that admit no PNEs}.
We investigate this behaviour in the context of random games with 
$N$ players, each choosing between two actions and receiving independently and identically distributed  payoffs. The payoff distribution is permitted to have atoms, allowing for the possibility of ties. The induced dynamics define a class of random walks on the $N$-dimensional hypercube, where the landscape includes two types of absorbing sets: PNEs and traps. Our analysis encompasses BRD, simple random walks, and a broad class of nearest-neighbour processes. We prove that, with high probability, these dynamics reach a PNE before encountering any trap. Our results complete and generalise the ones in \cite{amiet2019pure} by Amiet, Collevecchio, Scarsini and Zhong, using a different approach.
\end{abstract}

\maketitle

\tableofcontents

\section{Introduction.}
\subsection{Motivation}
Understanding the typical structure of games is fundamental to the development of efficient algorithms for computing solutions, particularly Mixed Nash Equilibria (MNEs). In this setting, a MNE refers to a profile of {randomised} strategies in which no individual player has an incentive to unilaterally deviate. A special subclass of MNEs is formed by the Pure Nash Equilibria (PNEs), which lie directly in the space of pure strategies rather than in distributions over them. While every finite game is guaranteed to admit at least one MNE, it may possess no PNEs at all.
The reader can think of the familiar  Rock-Paper-Scissors game. 
Although the definitions and properties of mixed strategies and mixed equilibria are clear, the concept of PNEs is ``more natural". It is worth mentioning that Osborne and Rubinstein discussed the interpretation of mixed equilibria in Section 3.2 of \cite{OR1994} with two different paragraphs individually signed by the two authors as they could not reach an agreement. Generally, PNEs exhibit a stronger epistemic foundation than mixed equilibria. As mentioned earlier, the main issue of PNEs is their existence.
In order to understand the typical behaviour of a game, it is natural to consider what happens in the randomised case, where the payoffs are random. We consider the simplest model in which players’ payoffs are {independently and identically distributed}. It is natural to address the following questions. 

 $\bullet$ \textit{How many PNEs does a random game with $N$ players, each having two actions, typically have?} Answers to this question were given in the literature as follows.
\begin{enumerate}
\item \cite{RinSca:GEB2000, McLennan, Sta:EL1999, Sta1995, Sta:MOR1996} focused on the case when the payoffs are atomless, i.e. the probability that two payoffs are equal, which we denote by $\alpha$, is equal to zero. In this case the number of PNEs, which we denote by $\# {\rm PNEs}$ converges in distribution to a Poisson with mean one.
\item \cite{amiet2019pure} studied the case where ties are allowed, i.e. $\alpha \in (0,1)$. In this case 
$$\frac{(\# {\rm PNEs}) - (1+\alpha)^N}{(1+\alpha)^{N/2}}$$ 
converges in distribution to a standard normal.
This is used to prove that $\# {\rm PNEs}$ is concentrated around $(1+\alpha)^N$. For example,  consider a win-lose game, where the payoffs are i.i.d. $\text{Bernoulli}(1/2)$, with 100 players. In this case $\alpha = 1/2$, and the number of PNEs is ``very close" to $(1.5)^{100}$. 
\end{enumerate}

$\bullet$ \textit{How to find PNEs?} 
A natural way is to devise iterative procedures that converge to a PNE. 
For instance, some adaptive procedures start from a strategy profile and allow a single player (picked \emph{at random}) to choose a different strategy.
In the case of Best Response Dynamics (BRD), introduced by Gilboa and Matsui \cite{GM1991, Matsui}, the selected player myopically switches to a {``best response" strategy given the strategies chosen by other players}. The new strategy profile corresponding to this best response gives the highest possible payoff to the player. 
If no such strategy exists, the player will not move, and a different player will be chosen at random. 
The process is either absorbed by a PNE, if one exists, or may get stuck in subgraphs called \lq traps' or \lq sink equilibria' (see, e.g., \cite{GoeMirVet:FOCS2005}).
The following results for random games with $N$ players, each having two actions, were proved in \cite{amiet2019pure}. \begin{itemize}
    \item[i.]  If $\alpha =0$, then BRD  converges to a PNE conditionally on the existence of the latter, which holds with probability $1- {\rm e}^{-1} + o(1)$.
   \item[ii.] If $0<\alpha < 2^{3/4} -1\approx 0.6818$, then the BRD converges to a PNE with high probability, as shown in \cite{amiet2019pure}, as there are many PNEs, as mentioned above, and there are no traps, i.e. subgraphs that prevent BRD from reaching PNEs, with very high probability.
\end{itemize}

For a formal definition of PNEs, traps and BRD see Section~\ref{se:number}.

\begin{figure}[ht]   
\centering 
\begin{tikzpicture}
  \begin{axis}
    [
    xlabel={$\alpha$},
    boxplot/draw direction=y,
    xtick={1,2,3,4,5},
    xticklabels={0.5, 0.6, 0.7, 0.8, 0.9 },
    ytick = {0, 50, 100, 150, 200, 250 },
    ymajorgrids,
    boxplot={variable width,
            box extend=0.3,
            draw position={0.8 + floor(\plotnumofactualtype/2) + 1/3*mod(\plotnumofactualtype,2)}
            }
    ]
    % Best response dynamics
    \addplot+[
    boxplot prepared={
      median=39,
      upper quartile=77,
      lower quartile=15,
      upper whisker=166,
      lower whisker=3,
      every median/.style={ultra thick,solid, draw=red},
      every whisker/.style={solid,draw=red},
      every box/.style={solid,draw=red},
    },
    fill=red!30,
    draw=red,
    ] coordinates {};
    % Simple random walks
    \addplot+[
    boxplot prepared={
      median= 60,
      upper quartile=118.5,
      lower quartile=23,
      upper whisker=238.5,
      lower whisker=3,
      every median/.style={ultra thick,solid, draw=blue},
      every whisker/.style={solid,draw=blue},
      every box/.style={solid,draw=blue},
    },
    fill=blue!30,
    draw=blue,
    ] coordinates {};
    % Add the second set of boxplots for alpha=0.6
    \addplot+[
    boxplot prepared={
      median=15,
      upper quartile=6,
      lower quartile=31,
      upper whisker=64,
      lower whisker=1,
      every median/.style={ultra thick,solid, draw=red},
      every whisker/.style={solid,draw=red},
      every box/.style={solid,draw=red},
    },
    fill=red!30,
    draw=red,
    ] coordinates {};
    \addplot+[
    boxplot prepared={
      median=21,
      upper quartile=8,
      lower quartile=40.5,
      upper whisker=88,
      lower whisker=1,
      every median/.style={ultra thick,solid, draw=blue},
      every whisker/.style={solid,draw=blue},
      every box/.style={solid,draw=blue},
    },
    fill=blue!30,
    draw=blue,
    ] coordinates {};
    % Add the third set of boxplots for alpha=0.7
    \addplot+[
    boxplot prepared={
      median=6,
      upper quartile=12,
      lower quartile=2,
      upper whisker=27.5,
      lower whisker=0,
      every median/.style={ultra thick,solid, draw=red},
      every whisker/.style={solid,draw=red},
      every box/.style={solid,draw=red},
    },
    fill=red!30,
    draw=red,
    ] coordinates {};
    \addplot+[
    boxplot prepared={
      median=8,
      upper quartile=17,
      lower quartile=3,
      upper whisker=36,
      lower whisker=0,
      every median/.style={ultra thick,solid, draw=blue},
      every whisker/.style={solid,draw=blue},
      every box/.style={solid,draw=blue},
    },
    fill=blue!30,
    draw=blue,
    ] coordinates {};
    % Add the fourth set of boxplots for alpha=0.8
    \addplot+[
    boxplot prepared={
      median=2,
      upper quartile=5,
      lower quartile=1,
      upper whisker=11,
      lower whisker=0,
      every median/.style={ultra thick,solid, draw=red},
      every whisker/.style={solid,draw=red},
      every box/.style={solid,draw=red},
    },
    fill=red!30,
    draw=red,
    ] coordinates {};
    \addplot+[
    boxplot prepared={
      median=3,
      upper quartile=6,
      lower quartile=1,
      upper whisker=12,
      lower whisker=0,
      every median/.style={ultra thick,solid, draw=blue},
      every whisker/.style={solid,draw=blue},
      every box/.style={solid,draw=blue},
    },
    fill=blue!30,
    draw=blue,
    ] coordinates {};
    % Add the fifth set of boxplots for alpha=0.9
    \addplot+[
    boxplot prepared={
      median=1,
      upper quartile=1.5,
      lower quartile=0,
      upper whisker=4,
      lower whisker=0,
      every median/.style={ultra thick,solid, draw=red},
      every whisker/.style={solid,draw=red},
      every box/.style={solid,draw=red},
    },
    fill=red!30,
    draw=red,
    ] coordinates {};
    \addplot+[
    boxplot prepared={
      median=1,
      upper quartile=2,
      lower quartile=0,
      upper whisker=5,
      lower whisker=0,
      every median/.style={ultra thick,solid, draw=blue},
      every whisker/.style={solid,draw=blue},
      every box/.style={solid,draw=blue},
    },
    fill=blue!30,
    draw=blue,
    ] coordinates {};
  \end{axis}
\end{tikzpicture}
\caption{Iterations needed for Best Response Dynamics (red) and Simple Random Walk (blue) to reach a PNE conditional on not hitting any trap for $N=15$ and $\alpha\in \{0.5, 0.6, 0.7, 0.8, 0.9 \}$, with 500 trials. The boxes and the whiskers are corresponding to (0.25, 0.75) and (0.05, 0.95) quantile intervals respectively.}\label{fig:BRD_step}
\end{figure}
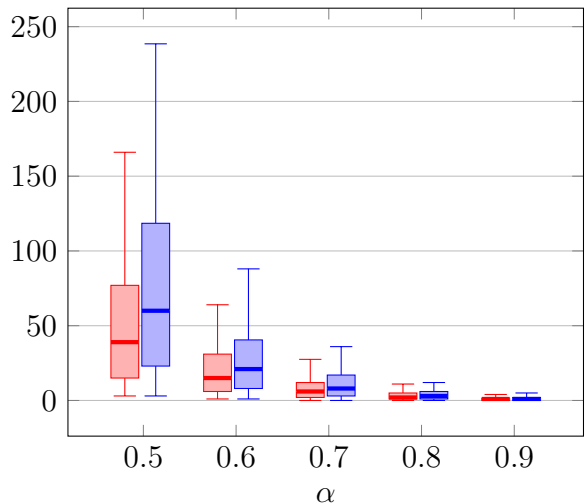

The box plot in \cref{fig:BRD_step} shows that as $ \alpha $ increases, the median number of iterations decreases, indicating faster convergence to a PNE, conditional on not hitting a trap. This phenomenon is driven by the expected number of PNEs, given by $ (1+\alpha)^N $, which grows exponentially with $\alpha$. The plot also suggests that BRD converges to a PNE typically faster than the simple random walk.

The goal of this paper is to fill the gap in \cite{amiet2019pure}, and study in detail the case where there is a competition between PNEs and traps. It is worth mentioning that when the tie probability $\alpha$ increases, the probability that BRD hits a trap also increases as there are more traps.~~In our paper, we prove that BRD still reaches a PNE before hitting any trap with high probability, for any $\alpha \in (0,1)$. Our result (see \cref{th:BRD} below) holds for a large class of nearest neighbour processes. { While finding a Nash equilibrium is hard in the worst case (PPAD-complete), random games represent a “typical” case. The result therefore highlights a separation between worst-case hardness and average-case tractability. It provides evidence that equilibria may be behaviourally plausible: typically, in large games, even uncoordinated agents following simple improvement rules can converge to stability.}

\subsection{Further background}
 The probability that a random two-person game contains a PNE was first studied for zero-sum games in \cite{Goldman} and for games with i.i.d. payoffs in \cite{Goldberg}. The result in \cite{Goldberg} was later extended to $N-$person games with i.i.d. payoffs in
\cite{Dresher} and \cite{Papavassilopoulos}. Moreover, it was shown that if the distribution of the payoffs does not admit atoms, then the number of PNEs is asymptotically Poisson(1) when the dimension of the game grows (see, for example,  \cite{Powers, Sta1995, RinSca:GEB2000}). 
Much more difficult to study is the case where the distribution admits atoms. This case was studied in \cite{amiet2019pure}, where the  authors obtained a central limit theorem for the number of PNEs, with a description of the rate of convergence.  Other models with a large number of PNEs related to vector payoffs and common payoffs were investigated in \cite{Sta:MSS1997} and \cite{Sta:EL1999}.
  It is worth mentioning that the framework of games with random payoffs was extended to graphical games in \cite{Daskalakis}, in which a graph is defined whose vertices represent the players of the game and edges correspond to
the strategic interaction between players. The payoff of  vertex $v$ depends only on the strategies of $v$ and its neighbours. The authors in \cite{Daskalakis} investigated how the structure of the graph affects the existence of PNEs and studied the convergences to Poisson(1) of the number of PNEs. The relation between the one-shot deviation property and the stochastic stability of random games having a PNE was established in \cite{newton2024conventions}. Recently, \cite{alon2021dominance} studied the connection between the dominance solvability of random games and the iterated elimination of pure-strategy strictly dominated actions required for the convergence to a PNE.  

As for the adaptive processes on random games, the behaviour of BRD and better response dynamics (bRD) on games with i.i.d. atomless payoffs,  two players, and a large number of actions was studied in \cite{Amiet2021WhenI}. In particular, the authors in \cite{Amiet2021WhenI} established that BRD does not converge to a PNE, whereas bRD does whenever at least one PNE exists. More precisely, BRD is absorbed in a trap of length $\sqrt{N}$, where $N$ is the number of actions, with high probability, for $N$ large enough. Additionally, the behaviour of BRD in random games with correlated payoffs  was studied in \cite{Mimun}. The authors in \cite{Mimun} described a phase transition in the convergence of BRD, in terms of the correlation parameter of the payoffs. It is worth mentioning that a connection between the structural properties of games and the connectivity properties of their best-response graphs has been recently studied in \cite{JSST2023}. In particular, the authors in \cite{JSST2023} have substantially extended the results in \cite{Amiet2021WhenI} to games with a large number of players, each having multiple actions and i.i.d. atomless payoffs. 

The case of games with $N$ players, each having two actions and i.i.d. payoffs with atoms was treated in \cite{amiet2019pure}, where the authors proved that BRD converges whenever $\alpha \in (0, 2^{3/4} -1)$ where $\alpha$ is the probability that two payoffs coincide.
In the present paper we show that BRD reaches a PNE before hitting any trap for all $\alpha \in (0,1)$ and for all large enough games, i.e., when $N \to \infty$. 
 Our result holds not only for BRD but also for a large class of processes.

\subsection{Basic concepts from Game Theory}\label{se:number}

A game $\Gamma$ is a triplet $\big(\mathcal{P},(S_{i})_{i\in \mathcal{P}},( g _{i})_{i\in\mathcal{P}} \big)$,
where $\mathcal{P}$ is the set of players, $S_{i}$ is the set of strategies of each player $i\in \mathcal{P}$ and  {$ g _{i}\colon \bigtimes_{j\in \mathcal{P}} S_j \to\R$ {is} the payoff function of player $i$}.

Let $N$ be a natural number which is not smaller than 2. Set $[N]:=\{1,2,\dots, N\}$.   We assume throughout this paper that $\Gamma=\Gamma_N$ is a game with $N$ players and each player has exactly two strategies. More specifically, we set $\mathcal {P}=[N]$ and $S_{i}=\{0,1\}$ for each $i\in [N]$. Set $\mathcal{V}_N := S_1\times S_2\times\cdots\times S_N = \{0,1\}^N$. For each $i\in [N]$, let $g_i^{(N)}:\vertices\to \R$ be the payoff function of player $i$. An element $u=(u_1,u_2,\cdots, u_N)\in \mathcal{V}_N$ is called a \textit{strategy profile}. For each $i \in  [N]$, denote by $u^{(i)}$ the strategy profile that differs from $u$ exactly in the action of the $i$-th player.
A pure Nash equilibrium in a normal-form game is a profile of strategies, one for each player, such that, given the choices of the other players, no player has an incentive to make a different choice. In other words, deviations from the equilibrium are not profitable for any player.

\begin{definition}\label{def.pne}
A strategy profile $u\in \mathcal{V}_N$ is called  Pure Nash Equilibrium (PNE) of the game $\Gamma_{N}$ if 
	$$ g_{i}^{(N)}(u)\ge g_{i}^{(N)}(u^{(i)}) \text{ for all  } i \in  [N].$$
\end{definition}

An example with two players is given in \cref{fi:game3}.  In this example, the strategy profile $(0,0)$ is a  PNE. We can identify the game's strategies as the vertices of a hypercube $\cube=\parens{\vertices,\edges}$ where an edge connects two vertices if they differ in exactly one coordinate.  Define the \textit{Hamming distance} $ d \colon \vertices\times \vertices \mapsto \{0, 1, \ldots N\}$ as
$$
 d(u,v)\coloneqq \card{ i\in [N]  \colon  u_{ i }\neq v_{ i }},
$$
 \begin{wraptable}{r}{9cm}\centering
	\begin{tabular}{|c|c|c|}
%		\hline
%		\multicolumn{3}{|c|}{Player 3 - Strategy  $0$}\\ 
		\hline
		\diagbox{P1}{P2} & $0$ & $1$\\ 
		\hline
		$0$ & $(0.322,0.412)$ & $(0.469,0.233)$ \\
		\hline
		$1$ & $(0.214,0.878)$ & $(0.202,0.311)$ \\
		\hline
%		\hline
%		\multicolumn{3}{|c|}{Player 3 - Strategy $1$}\\ 
%		\hline
%		\diagbox{Player 1}{Player 2} & $0$ & $1$\\  
%		\hline
%		$0$ & $(0.815,0.774,0.508)$ & $(0.292,0.684,0.126)$ \\ 
%		\hline
%		$1$ & $(0.209,0.659,0.4)$ & $(0.542,0.709,0.426)$ \\
%		\hline
	\end{tabular}
	\caption{Representation of a game $\Gamma_{2}$ on  $\braces{0,1}^{2}$.}
	\label{fi:game3}		
\end{wraptable}
 in which for a set $A$, we denote by $\card{A}$ the cardinality of $A$.
 We write  $u\sim v$ to denote that $u$ and $v$  are \emph{neighbours} in $\cube$, i.e., are connected by an edge which is equivalent to having Hamming distance one.  A subset of vertices $W\subseteq \vertices$ is \textit{connected} if for any distinct vertices $u, v\in W$ there is a path $(x_0,x_1,\dots, x_k)$ such that $x_0=u$, $x_k=v$ and $\{x_{i-1},x_{i}\}\in \mathcal{E}_N$ for all $ i\in [k]$. 
The vertex $(0,0, \ldots, 0) \in \vertices$ is denoted by $\zero$.

In order to identify PNEs, it is sufficient to know the ranking among the payoffs and not their specific values. This suggests that it is enough to describe the game using a partially oriented hypercube $\overrightarrow{H}(\Gamma_N)$, as follows: to each pair of neighbours  $u,u^{(i)}$ is associated either an oriented edge pointing from $u$ to $u^{(i)}$ if $g_i^{(N)}( u) < g_i^{(N)}(u^{(i)})$, denoted by $(u, u^{(i)})$,  or an unoriented edge if $g_i^{(N)}( u) = g_i^{(N)}(u ^{(i)})$. In this context,  a PNE is a vertex $u$ of $\overrightarrow{H}(\Gamma_N)$ which is incident only to edges that are either unoriented or pointing towards $u$.

For an oriented subgraph $\overrightarrow{{G}}$ of $\overrightarrow{H}(\Gamma_N)$, we denote by $\vertices(\overrightarrow{{G}})$ and $\oedges(\overrightarrow{{G}})$ the sets of vertices and oriented edges of $\overrightarrow{{G}}$ respectively. For a pair of distinct vertices $u, v\in \vertices(\overrightarrow{{G}})$, we say that $v$ is \textit{accessible} from $u$ in $\overrightarrow{{G}}$ if there is an oriented path $(x_0,x_1,\dots, x_k)$ such that $x_0=u$, $x_k=v$ and $(x_{i-1},x_{i})\in \oedges(\overrightarrow{{G}})$ for all $ i\in [k]$. We say $\overrightarrow{{G}}$
 is \emph{strongly connected} if every vertex is accessible from every other vertex in that graph. 

A \emph{trap} is a strongly connected induced oriented subgraph $\overrightarrow{{T}}$ of  $\overrightarrow{H}(\Gamma_N)$ with four or more vertices, such that, for all $ u \in \vertices(\overrightarrow{{T}})$ and all $v \notin \vertices(\overrightarrow{{T}})$, we have that $v$ is not accessible from $u$ in $\overrightarrow{H}(\Gamma_N)$. As $\cube$ is a bipartite graph, the smallest cycles of the graph are squares, i.e., cycles of length four. \cref{fig:4-cube} shows an example of a game with four players, in which there is a trap of size six. 

\begin{figure}[ht]\centering     
%\begin{subfigure}{0.4\textwidth}
\begin{tikzpicture}[scale=.6]
    	\scriptsize
    	\node [label={below left:$0$}] (0) at (0,0) {};
    	\node [label={below right:$1$}] (1) at (4,0) {};
    	\node [label={left:$2$}] (2) at (0,4) {};
    	\node [label={right:$3$}] (3) at (4,4) {};
    	\node [label={below:$4$}] (4) at (2.55,2.55) {};
        \node [label={left:$8$}] (8) at (-2.55,2.55) {};
        \node [label={below:$9$}] (9) at (1.275,2.55) {};
        \node [label={above left:$10$}] (10) at (1.275,6.55) {};
        \node [label={left:$11$}] (11) at (-2.55,6.55) {};
        \node [label={left:$12$}] (12) at (0,5.1)  {};
        \node [label={above left:$13$}] (13) at (0,9.1)  {};
        \node [label={above right:$14$}] (14) at (4,9.1)  {};
        \node [label={right:$15$}] (15) at (4,5.1)  {};
    	\node [label={below right:$5$}] (5) at (6.55,2.55) {};
    	\node [label={above right:$6$}] (6) at (2.55,6.55) {};
    	\node [label={above right:$7$}] (7) at (6.55,6.55) {};      
    	\filldraw (0,0) circle (3pt) [color=green];
    	\filldraw (4,0) circle (3pt) [color=black];
    	\filldraw  (0,4) circle (3pt) [color=black];
    	\filldraw (4,4) circle (3pt) [color=black];
    	\filldraw (2.55,2.55) circle (3pt) [color=red];
    	\filldraw(-2.55,2.55)circle (3pt)[color=black];
    	\filldraw (1.275,2.55)  circle (3pt) [color=black];
    	\filldraw (1.275,6.55) circle (3pt) [color=black];
    	\filldraw (-2.55,6.55) circle (3pt) [color=black];
    	\filldraw (0,5.1) circle (3pt) [color=black];
        \filldraw (0,9.1) circle (3pt) [color=black];
        \filldraw (4,9.1) circle (3pt) [color=black];
        \filldraw (4,5.1)  circle (3pt) [color=black];
    	\filldraw (6.55,2.55) circle (3pt) [color=black];
    	\filldraw (2.55,6.55) circle (3pt) [color=black];
    	\filldraw  (6.55,6.55) circle (3pt) [color=black];
    % 	\begin{scope}[thick,decoration={markings, mark=at position 0.33 with {\arrow{latex}}, mark=at position 0.66 with {\arrow{latex}}}] 
	\begin{scope}[thick,decoration={markings, mark=at position 0.66 with {\arrow{latex}}}] 
        	\draw[postaction={decorate}] (1) -- (0)[color=lightgray];
        	\draw[postaction={decorate}] (0) -- (2)[color=green];
        	\draw[postaction={decorate}] (0) -- (4)[color=green];
        	\draw[postaction={decorate}] (1) -- (3)[color=lightgray];
        	\draw[postaction={thick}] (1) -- (5)[color=lightgray];
        	\draw[postaction={decorate}] (2) -- (3)[color=green];
        	\draw[postaction={decorate}] (6) -- (2)[color=lightgray];
        	\draw[postaction={decorate}] (3) -- (7)[color=green];
        	\draw[postaction={decorate}] (5) -- (7)[color=red];
        	\draw[postaction={decorate}] (6) -- (4)[color=lightgray];
        	\draw[postaction={decorate}] (5) -- (7)[color=lightgray];
        	\draw[postaction={decorate}] (6) -- (7)[color=lightgray];
        	\draw[postaction={thick}] (0) -- (8)[color=lightgray];
        	\draw[postaction={decorate}] (1) -- (9)[color=lightgray];
        	\draw[postaction={decorate}] (8) -- (11)[color=blue];
        	\draw[postaction={decorate}] (10) -- (9)[color=blue];
        	\draw[postaction={decorate}] (2) -- (11)[color=green];
        	\draw[postaction={thick}] (3) -- (10)[color=lightgray];
        	\draw[postaction={decorate}] (11) -- (10)[color=blue];
        	\draw[postaction={decorate}] (9) -- (8)[color=blue];
        	\draw[postaction={decorate}] (12) -- (8)[color=lightgray];
        	\draw[postaction={decorate}] (12) -- (4)[color=lightgray];
        	\draw[postaction={decorate}] (13) -- (11)[color=blue];
        	\draw[postaction={decorate}] (13) -- (6)[color=lightgray];
        	\draw[postaction={decorate}] (12) -- (13)[color=lightgray];
        	\draw[postaction={decorate}] (14) -- (13)[color=blue];
        	\draw[postaction={decorate}] (7) -- (14)[color=green];
        	\draw[postaction={decorate}] (10) -- (14)[color=blue];
        	\draw[postaction={decorate}] (12) -- (15)[color=lightgray];
        	\draw[postaction={decorate}] (15) -- (14)[color=lightgray];
        	\draw[postaction={decorate}] (15) -- (5)[color=lightgray];
        	\draw[postaction={decorate}] (15) -- (9)[color=lightgray];
        	\draw[postaction={decorate}] (5) -- (4)[color=lightgray];
    	\end{scope}
    \end{tikzpicture} 
    \caption{Representation of $\Gamma_{4}$ on $\braces{0,1}^{4}$. Blue vertices form a trap, red vertices represent PNEs and green lines correspond 
to possible best-response paths.}\label{fig:4-cube}
\end{figure}

\begin{wrapfigure}{l}{8cm}
    \centering
	\begin{tikzpicture}[scale=.6]
		\centering
    	\scriptsize
    	\node [label={below left:$(0,0,0)$}] (0) at (0,0) {};
    	\node [label={below right:$(1,0,0)$}] (1) at (4,0) {};
    	\node [label={above left:$(0,1,0)$}] (2) at (0,4) {};
    	\node [label={right:$(1,1,0)$}] (3) at (4,4) {};
    	\node [label={left:$(0,0,1)$}] (4) at (1.8,1.8) {};
    	\node [label={below right:$(1,0,1)$}] (5) at (5.8,1.8) {};
    	\node [label={above left:$(0,1,1)$}] (6) at (1.8,5.8) {};
    	\node [label={above right:$(1,1,1)$}] (7) at (5.8,5.8) {};
    	\filldraw (0,0) circle (3pt) [color=black];
    	\filldraw (4,0) circle (3pt) [color=black];
    	\filldraw (0,4) circle (3pt) [color=black];
    	\filldraw (4,4) circle (3pt) [color=black];
    	\filldraw (1.8,1.8) circle (3pt) [color=red];
    	\filldraw (5.8,1.8) circle (3pt) [color=black];
    	\filldraw (1.8,5.8) circle (3pt) [color=black];
    	\filldraw (5.8,5.8) circle (3pt) [color=red];
    	\begin{scope}[thick,decoration={markings, mark=at position 0.66 with {\arrow{latex}}}] 
        	\draw[postaction={decorate}] (0) -- (1)[color=green];
        	\draw[postaction={decorate}] (2) -- (0)[color=lightgray];
        	\draw[postaction={decorate}] (0) -- (4)[color=green];
        	\draw[postaction={decorate}] (1) -- (3)[color=green];
        	\draw[postaction={thick}]    (1) -- (5)[color=lightgray];
        	\draw[postaction={decorate}] (3) -- (2)[color=lightgray];
        	\draw[postaction={decorate}] (6) -- (2)[color=lightgray];
        	\draw[postaction={decorate}] (3) -- (7)[color=green];
        	\draw[postaction={decorate}] (5) -- (4)[color=lightgray];
        	\draw[postaction={decorate}] (6) -- (4)[color=lightgray];
        	\draw[postaction={decorate}] (5) -- (7)[color=lightgray];
        	\draw[postaction={decorate}] (6) -- (7)[color=lightgray];
    	\end{scope}
	\end{tikzpicture}
	\caption{\label{fig:3-cube}Representation of $\Gamma_{3}$ on $\braces{0,1}^{3}$. Red vertices are PNEs and green lines represent possible best-response paths.}
\end{wrapfigure}

We define a \textit{Best Response Dynamic} (BRD) process as follows. It is a Markov chain on the oriented hypercube  $\overrightarrow{H}(\Gamma_N)$. It starts from ${\bf 0}$ and jumps to the nearest neighbours until reaching a PNE. At each step, the process chooses an edge uniformly at random among the ones incident to the present position and pointing out from it, then the process moves towards that direction. The process stops when it reaches a PNE. To see an example, the reader can refer to \cref{fig:3-cube}. The red dots are PNEs of the game. The process starts at ${\bf 0}$, and it follows one of the two blue edges, which share the same probability of being chosen. This is equivalent to picking a player randomly and asking if they want to change their strategy, upon which they agree only if they would receive a better payoff.

\subsection{Random games}
We aim to study random games with $N$ players, each having two actions with i.i.d. payoffs. As mentioned above, our goal is to prove that a large class of random walks in random media on the hypercube $\cube$, including BRD, eventually reaches a {PNE} before hitting any trap, with a probability quickly tending to 1 as the number of players goes to infinity. By appealing to the first Borel-Cantelli lemma, the convergence of these processes to a PNE holds for all large $N$. 

Consider a probability space $(\Omega,\mathcal{F},\Prob)$, on which the following sequence of random games is defined. For each integer $N\ge 2$, let {$\Gamma_{N}$} be a game with $N$ players, each with two strategies available, and random i.i.d.\ payoffs. Formally, the payoffs $(g_{i}^{(N)}(u))_{i\in [N], u\in\vertices}$ are i.i.d. copies of a random variable $Z$. 
Let 
$$\alpha\coloneqq\Prob\parens*{Z_{1}=Z_{2}}\quad\text{and}\quad\beta \coloneqq (1-\alpha)/2,$$
where $Z_{1}$ and $Z_{2}$ are i.i.d.\ copies of $Z$. 

{We stress that the random games $(\Gamma_N)_{N \geq 2}$ are not necessarily independent. 
More generally, we may consider any coupling between these games. For instance, let $(\Gamma_N)_{N \geq 2}$ be a sequence of random games with payoff functions defined as follows. Let $(Z_i^{x})_{i \in \mathbb{N},\; x \in \{0,1\}^{\mathbb{N}}}$ be i.i.d. copies of the random variable $Z$. For each $N \geq 2$, the payoffs in $\Gamma_N$ are given by $g_i^{(N)}(u) = Z_i^{\widetilde{u}}$ for $i \in [N]$ and $u \in \vertices$,
where $\widetilde{u} \in \{0,1\}^{\mathbb{N}}$ has first $N$ entries equal to $u$ 
and the rest zero.}

We associate the random game ${\Gamma_N}$ with a random partially oriented graph $\ocube = (\vertices, \oedges)$, as described in the previous section. Recall that $(u, u^{(i)}) \in \oedges$ if and only if $g_i^{(N)}(u) < g_i^{(N)}(u^{(i)})$. {Note that this property follows from the two-action-per-player assumption.} The orientation of each edge in $\oedges$, uniquely determined by a distinct pair of random variables $(g_{i}^{(N)}(u))_{i\in [N], u\in\vertices}$, is independent of the others. The set of oriented edges $\oedges$ can be viewed as a random medium induced by $\Gamma_N$ on the unoriented hypercube $\cube$.
%We provide a proof to the edge-independent property of $\ocube$ in the Appendix, \cref{le: edge_indep}.
We denote by $\mathbb{Z}_+$ and $\mathbb{N}=\mathbb{Z}_+\setminus\{0\}$ respectively the set of non-negative integers and the set of natural numbers. %{\color{red}(can we use $\mathbb{N}^{+}$ instead?)}
Let ${\bf X}^{(N)}=( X_t^{(N)})_{t\in \Z_+}$ be a discrete-time random walk in the random medium induced by $\Gamma_{N}$, in which its distribution is defined as follows. The process takes values in $\mathcal{V}_N$. For each fixed outcome $\overrightarrow{H}$ of $\ocube$, we denote by $P_{\overrightarrow{H}}$ the conditional probability measure corresponding to $\mathbb P$ on the event $\{\ocube=\overrightarrow{H}\}$. The distribution of ${\bf X}^{(N)}$ under $P_{\overrightarrow{H}}$ is called the {\it quenched law}.
Let $\mathfrak{H}_N$ be the set of all possible outcomes of $\ocube$.
The distribution of ${\bf X}^{(N)}$, which is called the {\it annealed law}, is given by 
\begin{equation}\label{annealed}
    \Prob(\mathbf X^{(N)}\in A)=\sum_{\overrightarrow{H}\in \mathfrak{H}_N} P_{\overrightarrow{H}}(\mathbf X^{(N)}\in A) \mathbb P(\ocube=\overrightarrow{H})
    \end{equation}
for each Borel measurable subset $A$ of $\mathcal{V}_N^{\mathbb \Z_+}$.

\subsection{Main result}
From now on, we assume the quenched law of ${\bf X}^{(N)}$ satisfies the following conditions.

\begin{assumption}\label{assump}
For any fixed outcome $\overrightarrow{H}$ of $\ocube$, we assume that under $P_{\overrightarrow{H}}$, the process ${\bf X}^{(N)}$ is a time-homogeneous Markov chain taking values on $\mathcal{V}_N$ such that $P_{\overrightarrow{H}}(X^{(N)}_0=\mathbf{0})=1$ and its transition probabilities satisfy the following conditions. If $x$ is a PNE in $\overrightarrow{H}$, then $P_{\overrightarrow{H}}(X^{(N)}_{t} = x\;|\; X^{(N)}_{t-1} = x)=1$. Otherwise, $P_{\overrightarrow{H}}(X_{t}^{(N)}\sim X^{(N)}_{t-1}|\; X^{(N)}_{t-1} = x)=1$ and for each $y\sim x$,
 \begin{align}\label{quenched}
     P_{\overrightarrow{H}}(X^{(N)}_{t} = y\;|\; X^{(N)}_{t-1} = x) \ge \kappa_1 N^{-\kappa_2} \delta_{(x, y) \in \overrightarrow{H}},
     \end{align}
 where $\delta$ is the usual Dirac mass, $\kappa_1$ and $\kappa_2$ are positive constants which do not depend on $N$, $x, y$ and $\overrightarrow{H}$. {Furthermore, we assume that for each $x\in\mathcal V_N$, the probabilities
 $P_{\overrightarrow{H}}(X^{(N)}_{t}=y\mid X^{(N)}_{t-1}=x)$, $y\sim x$, depend on $\overrightarrow{H}$ only through the orientations of the edges incident to $x$.}
 \end{assumption}

\begin{example}  It is clear that the simple random walk (SRW) on the unoriented hypercube $\cube$ satisfies  \cref{assump}. 
\end{example}
\begin{example}  The class of processes described in \cref{assump} includes the following random walk on the (unoriented) connected component of $\ocube$ containing $\mathbf{0}$. Let $\lambda_N\in (0,1]$ be a parameter such that $\lambda_N\ge \kappa_1 N^{-\kappa_2}$ for some positive constants $\kappa_1$ and $\kappa_2$. Assume that, under the quenched law, the transition probabilities are given by
    $$P_{\overrightarrow{H}}\bigl(X_t^{(N)}=y \,\big|\, X_{t-1}^{(N)}=x\bigr)
=
\frac{
\lambda_N \delta_{(x,y)\in\overrightarrow{H}}
+
(1-\lambda_N)\delta_{(y,x)\in\overrightarrow{H}}
}{
\lambda_N \deg_{\mathrm{out}}(x)+(1-\lambda_N)\deg_{\mathrm{in}}(x)
},$$
whenever $x$ is not a PNE, where $\deg_{\mathrm{out}}(x)$ and $\deg_{\mathrm{in}}(x)$ denote, respectively, the number of outward directed edges from $x$ and the number of inward directed edges to $x$. In other words, when the process is not at a PNE, it moves along an outward directed edge with probability $\lambda_N$ and along an inward directed edge with probability $1-\lambda_N$. In particular, when $\lambda_N=1$, the process reduces to BRD, that is, it follows only outward directed edges, which correspond to best responses. {When $\lambda_N\in(0,1)$, the process is no longer identical to BRD, since it is also allowed to move against the orientation. In that case, traps are not necessarily absorbing. However, in the regime where $\lambda_N$ is close to $1$ as $N\to\infty$, traps may act as metastable regions: once the process enters one, the time required to leave it can be very large. Therefore, entering a trap before reaching a PNE may still substantially delay the hitting time of a PNE.}
\end{example}

\begin{remark}
The preceding examples illustrate that, for the class of processes covered by
\cref{assump}, PNEs are absorbing by definition, whereas traps may not be absorbing. The aim of this
paper is to establish a general result for nearest-neighbour processes showing
that a PNE is reached before any trap is encountered; see \cref{th:BRD} below.
We believe that this result sheds light on the geometry of PNEs and traps.
In ongoing work, we investigate the time required for these processes to reach
a PNE \cite{ColleMax}.
\end{remark}

Define $\pne_{N}$ to be the set of PNEs in $\ocube$.{ In particular $\pne_{N}$ is a random subset of $\vertices$}.  Let 
$$\tau^{(N)} := \inf\{t \colon X_t^{(N)} \in  \pne_{N}\}$$ 
be the first hitting time of the process  $\mathbf{X}^{(N)} $ to a PNE. Let $$\xi^{(N)} := \inf\{t \colon  X_t^{(N)} \in  \vertices(\overrightarrow{\trap}) \text{ for some trap } \overrightarrow{\trap} \text{ of } \ocube\}$$ 
be the first time the process $\mathbf{X}^{(N)} $ hits a trap.
{Since PNEs and traps are random oriented subgraphs that depend on the random orientation of the hypercube, the random times $\tau^{(N)}$ and $\xi^{(N)}$ are not necessarily stopping times with respect to the natural filtration of the process. This makes it more complicated to determine whether the process $\mathbf{X}^{(N)}$ reaches a PNE before encountering a trap.}

We say that a sequence of events $(A_N)_N\subset \Omega$ holds 
\begin{itemize}
    \item \textit{with small probability} (w.s.p.) if $ \sum_{N=1}^{\infty} \Prob(A_N)<\infty;$
    \item \textit{with high probability} (w.h.p.) if the sequence of complements $(A^c_N)_N$ holds w.s.p.
\end{itemize}
{By the first Borel-Cantelli lemma, if the sequence $(A_N)_N$ holds w.h.p.,  then there exists a sufficiently large (random) $N_0$ such that $A_N$ occurs for all $N \ge N_0$.}

The following is our main result:

\vspace{2mm}

\begin{theorem}\label{th:BRD} Let $\alpha\in (0,1)$, i.e., $\beta\in (0, 1/2)$. {Assume that for each $N$, the quenched law of the process $\mathbf{X}^{(N)}$ satisfies \cref{assump}.} Then the sequence of events $({\tau^{(N)}<\xi^{(N)}})_{N\ge 2}$ holds w.h.p. Consequently, almost surely there exists a random $N_0$ such that for every $N\ge N_0$, the process $\mathbf{X}^{(N)}$ hits a PNE before hitting any trap.
\end{theorem}

% For every subgraph $C$ of $\ocube$ we denote by $\orient{\mathcal{E}(C)}$ its oriented edge set while $\mathcal{E}(C)$ its unoriented edge set.
%The results in \cref{th:BRD}  are generalised in \cref{th:BRDconv} below, where information about the rate of convergence is provided.
%\begin{mdframed}
%[style=MyFrame1]
%\begin{theorem}
%\label{th:BRDconv} 
%Fix $c, \alpha\in (0, 1]$.  Assume  $\mathbf{X}^{(N)} \in \mathcal{G}_N(c)$, for all $N\in \N$.  For any $\epsilon>0$, we have
%\begin{equation}
%\label{eq:rate of converg}
%\sum_{N=1}^{\infty} \Prob\left(\tau^{(N)} \le (1- \epsilon) \left(\frac 2{1+\alpha}\right)^N \right)<\infty.
%\end{equation}
%\end{theorem} 
%
 
\begin{remark} \cref{th:BRD} generalises the result in \cite{amiet2019pure}, where the authors proved that, with high probability, BRD converges to a PNE for  $0<\alpha<2^{3/4}-1\approx 0.68$. %, by using different tools. In fact, when $\alpha \in (0,2^{3/4}-1)$, it is very likely to have at least one {PNE} and no traps, which in turn implies that with high probability, the best response dynamics converges to a PNE. 
In contrast, it was shown in \cite{heinrich2021best} that if at each step the players are not chosen at random, but they follow a fixed deterministic cyclic order, then the process, called \lq clockwork' BRD, does not converge to a PNE for either a large number of players or a large number of strategies. Moreover, \cite{Amiet2021WhenI} studies a two-player game with {$M$ actions} where the payoffs are i.i.d. continuous random variables, so that 
the game is represented by an $M \times M$ payoff matrix. In contrast to the two-action case,  when $M$ is large, the 
probability that BRD converges to a PNE tends to zero. Intuitively, with many actions, 
each player’s best response is unique with high probability, so BRD alternates 
deterministically between players and is very likely to revisit a previously seen row or 
column before reaching a PNE, thereby becoming trapped in a cycle,  with high probability.

In our model, when the tie probability is large, both the expected number of traps and the expected number of PNEs grow exponentially as the number of players (or equivalently, the dimension of the hypercube) increases.
%both the expected number of traps and the expected number of PNEs growth exponentially as the number of players (or equivalently, the dimension of the hypercube) increases. Using a coupling with bond percolation on hypercubes, 
We will show that traps and PNEs tend to cluster together geometrically. This geometric clustering phenomenon reveals that traps and PNEs typically coexist within the same local regions of the hypercube, providing a natural baseline for the analysis of more structured random games, including those with correlated payoffs or network constraints. For instance, the models studied in \cite{heinrich2021best} and \cite{Mimun} rely on strictly continuous payoff distributions (which exclude ties), whereas our framework is significantly more general, explicitly allowing payoff ties with positive probability.  Moreover, this underscores the complexity of finding equilibria in large random games.
\end{remark}

\subsection{Strategies of the proof}

In \cref{subse: dual_ocube}, we couple the random oriented hypercube $\ocube$ with an i.i.d. bond percolation on the hypercube $\mathcal{H}_N$, where the geometry of the latter is well-described in \cite{McDScoWhi:2020}.
{It is worth noting that the coupling with bond percolation was previously introduced in \cite{amiet2019pure}. However, that work considered only the restricted regime $\alpha \in (0, 2^{3/4} - 1)$, in which, with high probability, at least one PNE exists and no traps are present. Consequently, the best-response dynamics converges to a PNE with high probability in this parameter range. In contrast, for larger values of $\alpha$, the hypercube contains a huge number of both traps and PNEs, which are roughly uniformly distributed at the global scale.
The present paper completes the picture by analysing the full interval $\alpha \in (0,1)$. Achieving this requires substantially new techniques, as a full understanding of the geometric arrangement of both PNEs and traps inside the hypercube becomes essential. In fact, by combining the percolation coupling with the results of \cite{McDScoWhi:2020}, we will show that PNEs and traps form geometrically connected clusters. Heuristically, the process $\mathbf{X}^{(N)}$ explores these clusters one by one until it reaches a PNE or trap. Even though there are many traps, the PNEs outnumber them, which in turn implies that in each cluster, the process is likely to hit a PNE rather than a trap with high probability.}

Our strategy for the proof of \cref{th:BRD} is (presented in \cref{proof of brd}) can be roughly described as follows. We cover the set of absorbing subsets such that each group of traps and PNEs which are close to each other within a certain distance will be included in a cover set called cluster. We show in \cref{lem.cluster} that the process $\mathbf{X}^{(N)}$
can visit a large number of clusters satisfying a `regular' property. We next demonstrate in \cref{le:spe1} that with high probability these clusters contain PNEs only. We then establish in \cref{le:spe2} a lower bound for the probability that upon a visit to the boundary of a regular cluster, the process will reach a PNE inside this cluster. The reciprocal of this lower bound is polynomial in $N$, but much smaller than the number of regular clusters. Combining the above-mentioned results, we deduce that with high probability, the process eventually hits a PNE before hitting any trap.

\section{Randomly oriented hypercube and bond percolation}\label{subse: dual_ocube}

In this section, we explore the structure of the randomly oriented hypercube associated with the random game and establish a coupling with an independent and identically distributed (i.i.d.) bond percolation model on the hypercube. This coupling allows us to leverage known results from percolation theory to describe the geometry of the oriented graph, particularly the distribution of PNEs and traps. Let 
$$\olargo=\{v\in \vertices:  \text{ there exists an oriented path in $\ocube$ starting from $\mathbf 0$ to $v$}\}\cup\{\mathbf{0}\},$$
and $$\rlargo=\{v\in \vertices:  \text{ there exists an oriented path in $\ocube$ starting from $v$ to $\mathbf 0$}\}\cup\{\mathbf{0}\}.$$
{Notice that $\olargo$ and $\rlargo$ have the same distribution. To gain insight into the structure of $\olargo$ and $\rlargo$, we couple $\rlargo$ with an i.i.d. bond percolation on the underlying unoriented hypercube $ \cube = (\vertices, \edges)$.} %This coupling was previously mentioned in  \cite{amiet2019pure}.

\vspace{.2cm}

\begin{definition}\label{def.perc}
Fix a  parameter $\beta\in(0,1)$.  Each edge in $\cube$ is \emph{open} with probability $\beta$ independently of other edges, and is \emph{closed} otherwise. The random subgraph with vertex set $\vertices$, obtained by deleting the closed edges, is called \textit{bond percolation} with parameter $\beta$. 
\end{definition} 

\vspace{.2cm}

Recall that we assume $\beta=(1-\alpha)/2\in (0,1/2]$ as $\alpha\in [0,1)$. Define the process 
	\begin{equation*}\label{eq:Pk+1}
	\mathcal{Q}_{1}=\braces{\mathbf{0}}\text{ and, for each }k\in\N,\
	\mathcal{Q}_{k+1}=\mathcal{Q}_{k}\cup\big\{ v\in \mathcal{Q}_{k}^c : \exists u\in \mathcal{Q}_{k} \text{ such that } (v,u)\in\oedges \big\} .
	\end{equation*}
Start with a bond percolation on $\cube$ with parameter $\beta$ that is independent of the random variables $(g_{i}^{(N)}(u))_{i\in [N], u\in\vertices}$. 
Call the resulting graph $\mathcal{B}_{1}$.	For every $k\in \N$ we will update $\mathcal{B}_{k}$ by changing the status of some edges at each stage,  in such a way that  $\mathcal{B}_{k+1}$ is still a bond percolation with parameter $\beta$. For each edge $e \in \edges$, we define $\mathcal{B}_{k}[e]$ the status (open or closed) of edge $e$ in $\mathcal{B}_{k}$. 
	We obtain  $\mathcal{B}_{k+1}$ from  $\mathcal{B}_{k}$,  by updating \emph{all and only} the edges  in  $\edges$ that connect an element of $\mathcal{Q}_{k}$ to an element of $\mathcal{Q}_{k}^c$. 
	More precisely, for any $u \in \mathcal{Q}_{k}$ and any  $w \in \mathcal{Q}_{k}^c$, with $u \sim w$, set
	\begin{equation*}\label{eq:Bk+1}
	\mathcal{B}_{k+1}[u, w] = 
	\begin{cases}
	\text{open,} & \text{if } (w, u)\in \oedges,\\
	\text{closed,} & \text{otherwise};
	\end{cases}
	\end{equation*}
	for all other edges $e \in \edges$, we have $\mathcal{B}_{k+1}[e] = \mathcal{B}_{k}[e]$.
 
Set $\mathcal{B}^{\beta}_N:=\mathcal{B}_{2^N}$. We note that $\mathcal{B}^{\beta}_{N}$ is not necessarily connected. Let $\mathcal{L}^{\beta,\mathbf{0}}_N$ be the set of vertices in the connected component of the percolation graph $\mathcal{B}^{\beta}_N$ that contains $\mathbf{0}$. By construction, we notice that $\mathcal{Q}_{2^N}$ is exactly equal to $\mathcal{L}^{\beta,\mathbf{0}}_N$ and that \begin{align}\label{o.comp}\mathcal{L}^{\beta,\mathbf{0}}_N = \overleftarrow{\mathcal{L}}^{\beta,\mathbf{0}}_N.\end{align} 

It is worth mentioning that the above-mentioned coupling with percolation was previously introduced in \cite{amiet2019pure}.
 The coupling allows us to give a detailed description of the geometry of PNEs. Denote by $\ularg$ the set of vertices in the largest connected component of the percolation $\mathcal{B}^{\beta}_{N}$. {The following result strengthens Proposition~2 in \cite{amiet2019pure}, stating that $\rlargo$ coincides with $\ularg$ with high probability.}

\vspace{.2cm}

\begin{lemma}
\label{le:finally}  
Assume that $\alpha \in [0,1)$, i.e. $\beta \in (0, 1/2]$. There exists a constant $\theta\in(0,1)$ such that for large enough $N$, 
$$\Prob\big(\ularg \neq \rlargo\big)\le \theta^{N}.$$ 
\end{lemma}

\vspace{.2cm}
\begin{proof}
Recall from \eqref{o.comp} that $\mathcal{L}^{\beta,\mathbf{0}}_N = \overleftarrow{\mathcal{L}}^{\beta,\mathbf{0}}_N$. Hence, it is sufficient to show that there exists a constant $\theta\in(0,1)$ such that for large enough $N$,
\begin{align}\label{iq1}
\Prob\big(\ularg \neq \mathcal{L}^{\beta,\mathbf{0}}_N\big)\le \theta^{N}.
\end{align}

We first consider the case $\beta\in (0,1/2)$. Let $\mathcal{M}^{\beta}_N:=\vertices\setminus \ularg$,  \textit{the fragment of the percolation} $\mathcal{B}^{\beta}_{N}$. Fix $\epsilon>0$ and set $$\zeta_N:=\Expect\left[\card{\mathcal{M}^{\beta}_N}\right]+\epsilon \sqrt{N(2(1-\beta))^N}.$$
As in equation (37) in the proof of \cite[Proposition~2]{amiet2019pure}, we have
\begin{align}\label{iq2}
    \Prob\big(\ularg \neq \mathcal{L}^{\beta,\mathbf{0}}_N \big)\le \frac{\zeta_N}{2^N}+ \Prob\left(\card{\mathcal{M}^{\beta}_N}\ge \zeta_N\right).
\end{align}
Furthermore, in view of Theorem 3 and inequality (2) in \cite{McDScoWhi:2020}, we have that for $\beta\in (0,1/2),$
\begin{align} \label{iq3}&\Expect\left[\card{\mathcal{M}^{\beta}_N}\right]=(2(1-\beta))^{N}\big(1+O(N(1-\beta)^{N})\big), \end{align}
and there exists a constant $\theta_1\in(0,1)$ depending only on $\epsilon$ and $\beta$ such that for large enough $N$,
\begin{align}
\label{iq4}  &\Prob\left(\card{\mathcal{M}^{\beta}_N}\ge \zeta_N\right)   \le \theta_1^N,
\end{align}
Combining \eqref{iq2} with \eqref{iq3} and \eqref{iq4}, we infer that \eqref{iq1} holds for all $\beta\in (0,1/2)$.

% Similarly as in the proof of Proposition~2 in \cite{amiet2019pure}, the case $\beta=1/2$ is proved by using the fact that $\mathcal{L}_N^{\beta}$ is monotone in $\beta$.

Finally, we consider the case $\beta=1/2$.  As in the proof of  \cite[Proposition~2]{amiet2019pure}, the fact that  $\mathcal{L}_N^{\beta'}\subseteq \mathcal{L}_N^\beta$ for $\beta'<\beta$ implies
$$
\Prob\bigl(\mathcal{L}_N^{1/2} \neq \mathcal{L}_N^{1/2,\mathbf{0}}\bigr)
\le
\Prob\bigl(\mathcal{L}_N^{\beta'} \neq \mathcal{L}_N^{\beta',\mathbf{0}}\bigr) \quad \text{for any $\beta'\in (0,1/2)$}.
$$
Therefore \eqref{iq1} also holds in this case.
\end{proof}

For $r\in \R$ and $v\in \mathcal{V}_N$, set $\text{Ball}_r(v)=\{u\in \mathcal{V}_N : d(u,v)\le r\},$
where we recall that $d$ is the {Hamming distance} on $\mathcal{H}_N$. Let
$$m_{\beta}\coloneqq\floor*{(-\log_2(1-\beta))^{-1}}.$$
In virtue of \protect{\cite[Theorem~2(a)]{McDScoWhi:2020}}, we note that, w.h.p.,  $m_{\beta}$ is an upper bound for the size of the smaller connected components in the bond percolation $\mathcal{B}_N^{\beta}$. {The next lemma states that, w.h.p., 
for every vertex $v$, the ball of radius $\delta N$ (for some $\delta > 0$) centred at $v$ contains at most $m_\beta$ vertices outside of $\olargo$ or $\rlargo$.}
%there is no ball of radius $\delta N$ (for some $\delta > 0$) containing strictly more than $m_{\beta}$ vertices outside $\olargo$ or $\rlargo$.

\vspace{0.2cm}

\begin{lemma}
\label{le:McD2}
Fix $\beta \in (0, 1/2]$.  
There exist $\delta>0$ and $\theta\in (0,1)$ such that for large enough $N$,
\begin{align}
\label{smallball1} 
     \mathbb{P}\left(\exists\; v \colon \card{\ball_{\delta N}(v)\setminus\rlargo} > m_{\beta} \right) &\le \theta^N\quad \text{and}\\
     \label{smallball2}
    \mathbb{P}\left(\exists\; v \colon \card{\ball_{\delta N}(v)\setminus\olargo}  > m_{\beta}\right) &\le \theta^N.
\end{align}
\end{lemma}	

\begin{proof}
   By the proof of Theorem~2(a) in \cite{McDScoWhi:2020}, there exist $\delta>0$ and $\theta_1\in (0,1)$ such that for large enough $N$,
$$\mathbb{P}\left(\exists\; v \colon \card{\ball_{\delta N}(v)\setminus\mathcal{L}^{\beta}_N} > m_{\beta}\right) \le \theta_1^N.$$ 
By \cref{le:finally}, there exists $\theta_2\in (0,1)$ such that for large enough $N$,
$$\Prob\big(\ularg \neq \rlargo\big)\le \theta_2^{N}.$$ 
Hence
$$\mathbb{P}\left(\exists\; v \colon \card{{{\ball}_{\delta N}}(v)\setminus\rlargo} > m_{\beta}\right)\le \theta_1^N+\theta_2^N.$$
This immediately implies \eqref{smallball1}. As $\olargo$ has the same distribution as $\rlargo$, we obtain \eqref{smallball2}. 
\end{proof}

Define
\begin{align}\nonumber
	\Abso_{N}&:=\big\{\mathcal{U}\subseteq\vertices\colon \text{either }\mathcal{U}=\{v\}\text{ for some }v\in\pne_{N} \text{or } \mathcal{U}=\vertices(\overrightarrow{\trap})
\\\nonumber &\quad \text{for some trap }\overrightarrow{\trap} \text{ of }\ocube\big\},\\ \label{Abso.star}
 \Abso^*_N&\coloneqq \big\{\mathcal{U}\subseteq\vertices\colon \mathcal{U}\in \Abso_{N}\text{ and }\card{\mathcal{U}}\le m_{\beta}\big\}
 \quad\text{and}\\ \nonumber
 {\sf Trap}_N^*
&:=
\Abso_N^*\setminus \{\{v\}:v\in\pne_N\},
\end{align}
where the elements of $\Abso_N$ are referred to as the \textit{absorbing sets}. Accordingly, $\Abso_N^*$ denotes the collection of absorbing sets of size at most $m_\beta$, and ${\sf Trap}_N^*$ denotes the subcollection of non-singleton elements, that is, the traps of size  at least four and at most $m_\beta$.
 The next result shows that, with high probability, each element of $\Abso_N$ contains no vertex belonging to $\rlargo$.

\begin{lemma}\label{le: GinFrag} Assume that $\alpha\in (0,1)$. There exists a constant $\theta\in (0,1)$ such that for large enough $N$,
$$ \Prob\Big(  \exists\ \mathcal{U} \in \Abso_{N} \colon 
  \mathcal{U} \cap  \rlargo\neq \emptyset \Big)\le \theta^N.$$
\end{lemma}

\vspace{.2cm}
\begin{proof} 
Our  first task is to  prove 
\begin{equation}\label{eq:inclus}
\{\exists\mathcal{U} \in \Abso_{N} : \mathcal{U}\cap \rlargo\neq \emptyset \}\subseteq \{ \olargo=\mathcal{U}, \text{ for some }\mathcal{U} \in \Abso_{N}\}.
\end{equation}
Assume that there is an element $\mathcal{U} \in \Abso_{N}$ such that $\mathcal{U}\cap \rlargo\neq \emptyset$. Each edge in $\ocube$ connecting a vertex of $\mathcal{U}$ to one in $\vertices \setminus \mathcal{U}$ is either unoriented or pointing towards $\mathcal{U}$. On the other hand, as $\mathcal{U}\cap \rlargo\neq \emptyset$, the vertex $\mathbf{0}$ is accessible in $\ocube$ from a vertex in $\mathcal{U}$. If $\mathbf{0}\notin \mathcal{U}$, then there are vertices $u'\in \vertices\setminus\mathcal{U}$ and $u\in \mathcal{U}$ such that $u$ and $u'$ are connected by an oriented edge pointing outwards $u$, which contradicts $\mathcal{U} \in \Abso_N$. Therefore, we must have $\mathbf{0} \in \mathcal{U}$, yielding that $\olargo=\mathcal{U}$. Hence \eqref{eq:inclus} is verified.

To complete the proof, it is sufficient to show that there exists a constant $\theta\in (0,1)$ such that for large enough $N$,
\begin{equation}\label{pro: frag_subset}
\Prob \Big( \olargo=\mathcal{U}, \text{ for some }\mathcal{U} \in \Abso_{N}\Big) \le \theta^N.
\end{equation}

Fix a small $\epsilon\in (0,1)$. In virtue of \cite[Theorem~2]{amiet2019pure}, we have that there exists a constant $\theta_1\in (0,1)$ depending only on $\epsilon$ and $\beta$ such that for large enough $N$,
\begin{align}\label{nPNE}
    \Prob\Big(\card{\pne_N}\le (1-\epsilon)(1+\alpha)^N\Big)\le \theta_1^N.
\end{align}

A PNE is said to be `good' if it is incident to at least $m_\beta +1$ incoming oriented edges in $\ocube$. Denote by $\pne_N^{*}$ the set of good PNEs. We argue next that there exists $\theta_2 \in (0,1)$ such that for large enough $N$,
\begin{equation}\label{eq:crbound}
\Prob(v\notin \pne_N^{*} \mid v\in\pne_N) \le \theta_2^N.
\end{equation}
In fact, if $v$ is a PNE, then each edge incident to $v$ is either oriented towards $v$ or unoriented. Conditionally on $v$ being a PNE, the status of each edge incident to $v$ is independent of the others, and it has probability $\beta/(1-\beta)$ of being oriented towards $v$. Hence \eqref{eq:crbound} is a direct consequence of the fact that $\Prob(v\notin \pne_N^{*} \mid v\in\pne_N)$ is equal to the probability that a binomial random variable with parameters ($N$, $\beta/(1-\beta)$) is equal or less than $m_\beta$.\\
We next prove that  
\begin{align}\label{good.pne}
    \Prob( \pne_N^*=\emptyset, \card{\pne_N}\ge (1-\epsilon)(1+\alpha)^N )\le (\theta_2^N)^{\floor{\frac{(1-\epsilon)}3(1+\alpha)^N}}
\end{align}
    for large enough $N$.  In order to show \eqref{good.pne}, we first show that on the event $\{\card{\pne_N}\ge (1-\epsilon)(1+\alpha)^N\}$, there exists $\frac{1-\epsilon}3 (1+\alpha)^N$ PNEs with mutual distance at least two. Note that the hypercube $\cube$ is a bipartite graph. We divide the vertex set $\vertices$ into two disjoint groups, `even' and `odd'. A vertex is `even' if its Hamming distance from $\zero$ is an even integer, and it is `odd', otherwise. On the event that there are at least $(1-\epsilon)(1+\alpha)^N$ PNEs, then at least  $\frac{(1-\epsilon)}3(1+\alpha)^N$ of them are all even, or all odd. Without loss of generality, we work on the event that there exist at least $\frac{(1-\epsilon)}3(1+\alpha)^N$ even PNEs and  let the reader repeat the argument for the odd case. Each pair of these PNEs have distance at least two. It means that they are not incident to common edges. Therefore, the events that these PNEs are good are independent. Hence conditional on the event $\{\card{\pne_N}\ge (1-\epsilon)(1+\alpha)^N\}$, the number of good PNEs is bounded below by a binomial random variable with parameters ($\floor{\frac{(1-\epsilon)}3(1+\alpha)^N}$,  $1-\theta_2^N$).
The probability that such a binomial variable is equal to zero can be quantified as 
$(\theta_2^N)^{\floor{\frac{(1-\epsilon)}3(1+\alpha)^N}}$. Hence \eqref{good.pne} is verified.

Next, we work on the event that there exists at least a good PNE, 
say $v$, and $\olargo=\mathcal{U}$ for some element ${\mathcal U}\in \Abso_{N}$. As $v$ is a good PNE, there exist $m_{\beta}+1$ neighbour vertices $(v_i)_{1\le i\le m_{\beta}+1}$ which are connected to $v$ by oriented edges pointing towards $v$. If $\olargo=\mathcal{U}\neq\{v\}$ 
 we must have that $v\in \vertices\setminus\olargo$. In that scenario, we have $(v_i)_{1\le i\le m_{\beta}+1}\subset  \vertices\setminus\olargo$, otherwise there is an oriented path from $\mathbf{0}$ to $v$, which leads to a contradiction that $v\in\olargo$. If $\olargo=\mathcal{U}=\{v\}$, it is obvious that $(v_i)_{1\le i\le m_{\beta}+1}\subset  \vertices\setminus\olargo$. As a result, we would have that $\card{{{\ball}_{{\delta N}}}(v)\setminus\olargo} > m_{\beta}$ for any $\delta>0$ and sufficient large $N$. Combining this fact with \cref{le:McD2}, we have that there exists $\theta_3\in (0,1)$ such that
\begin{align}\label{p.abs}
  \mathbb{P}\Big(  \olargo=\mathcal{U}, \text{ for some }\mathcal{U} \in \Abso_{N},\  \pne_N^{*}\neq\emptyset\Big)\le \theta_3^N,
\end{align}
for large enough $N$. Combining \eqref{nPNE}, \eqref{good.pne} and \eqref{p.abs} with  a union bound, we obtain that
$$\mathbb{P}\Big(  \olargo=\mathcal{U}, \text{ for some }\mathcal{U} \in \Abso_{N}\Big)\le \theta_1^N+(\theta_2^N)^{\floor{\frac{(1-\epsilon)}3(1+\alpha)^N}}+\theta_3^N.$$
for large enough $N$.
Hence \eqref{pro: frag_subset} is verified and this completes the proof.
\end{proof}

The following result is a corollary of \cref{le:McD2} and \cref{le: GinFrag}, which asserts that, with high probability, every element of $\Abso_{N}$ is either a PNE or the vertex set of a trap containing at most $m_{\beta}$ vertices.

\begin{lemma}\label{co:impoimpo}
Fix $\alpha\in (0,1)$. There exists a constant $\theta\in (0,1)$ such that for large enough $N$,
$$
\Prob\Big(\Abso_{N} \neq  \Abso^*_N\Big)\le\theta^N.
$$
\end{lemma} 

\begin{proof}
   {Notice that if $\Abso_{N} \neq \Abso^*_N$, then there
is an element $\mathcal{U}\in \Abso_{N}$ such that $\mathcal{U}$ is the vertex set of a trap with more than $m_{\beta}$ vertices, and one of the following two cases holds. Either $\mathcal{U}\cap \rlargo\neq \emptyset$ or $\mathcal{U}\cap \rlargo= \emptyset$. The latter case implies that for some $\delta>0$ and large enough $N$, there exists a vertex $v$ such that the ball of radius $\delta N$ centred at $v$ contains more than $m_{\beta}$ vertices.
Using a union bound together with \cref{le:McD2} and \cref{le: GinFrag}, we thus have that for all sufficiently large $N$,
\begin{align*}
    \Prob\Big(\Abso_{N} \neq \Abso^*_N\Big) 
    &\le \Prob\Big( \exists\ \mathcal{U} \in \Abso_{N} \colon
  \mathcal{U} \cap \rlargo\neq \emptyset \Big)\\
  &+\Prob\left(\exists\; v \colon \card{\ball_{\delta N}(v)\setminus\rlargo} > m_{\beta} \right)\\
  &\le \theta^N,
\end{align*}
where $\theta$ is some constant in $(0,1)$. This completes the proof.}
\end{proof}

\section{Proof of \texorpdfstring{\cref{th:BRD}}{Theorem 2.3}}\label{proof of brd}

\subsection{Strong construction of \texorpdfstring{$\mathbf{X}^{(N)}$}{X}}
Without loss of generality, we use the following strong construction of the process $\mathbf X^{(N)}$ throughout this section.

Let $\mathsf{vec}(\edges)=\{(u,v), (v,u) : \{u,v\}\in \edges \}$ denote the set of all possible oriented edges. Let $(\chi^{\vec{e}}_{k})_{k\in\Z_+, \vec{e}\in \mathsf{vec}(\edges)}$ be a collection of independent exponential random variables with parameter 1, which are  independent of $(g_{i}^{(N)}(u))_{i\in [N], u\in\vertices}$. For a fixed outcome $\overrightarrow{H}$ of $\ocube$,
let $(p_{\overrightarrow{H}}(x,y))_{(x,y)\in \mathsf{vec}(\edges)}$ be a transition probability matrix such that if $x$ is not a PNE in $\overrightarrow{H}$ then 
$$p_{\overrightarrow{H}}(x,y)\ge \kappa_1{N}^{-\kappa_2} \delta_{(x,y)\in \overrightarrow{H}},$$
otherwise $p_{\overrightarrow{H}}(x,x)=1$.
We assume that for each $x\in\vertices$, the probabilities $p_{\overrightarrow H}(x,y)$, $y\sim x$, depend on $\overrightarrow H$ only through the orientations of the edges incident to $x$.

Set the initial state $ X_0^{(N)} = \mathbf{0}$. Assume that $\ocube = \overrightarrow{H}$ and $X^{(N)}_{t}=u$ for some $t\in\N$. If $u$ is a PNE, we set $X_{t+1}^{(N)}=u$. Otherwise, for each neighbour $ v\sim u$, let $k^{(u,v)}_t$ be the number of crosses from $u$ to $v$ up to time $t$. We define the next state as follows. {We use the convention that $a/0=+\infty$ for $a>0$.} If there are two or more vertices $v\in \vertices$ achieving the same minimum of $\chi_{k^{(u,v)}_t}^{{(u,v)}}/p_{\overrightarrow{H}}(u, v)$ then we say that the construction fails and the process remains at $u$ forever. 
Otherwise, $$v_* := \arg\min_{v \in \vertices} \frac{\chi_{k^{(u,v)}_t}^{{(u,v)}}}{p_{\overrightarrow{H}}(u, v)}$$
is uniquely defined. In this case, we set the next state $X_{t+1}^{(N)} = v_*$.
Notice that the construction fails with probability 0 since the random variables $\chi_{k^{(u,v)}_t}^{{(u,v)}}/p_{\overrightarrow{H}}(u, v)$, with $v\sim u$ such that $p_{\overrightarrow{H}}(u, v)>0$,  are independent and have continuous distributions.
%where we use the convention that $ \chi_t^v / p_{\overrightarrow{H}}(u, v) = \infty $ if $ p_{\overrightarrow{H}}(u, v) = 0$.
%To ensure a well-defined construction, we specify a tie-breaking rule: order the vertices $ \vertices = \{0, 1\}^N $ \textit{lexicographically}, where $ u = (u_1, u_2, \ldots, u_N) $ precedes $ v = (v_1, v_2, \ldots, v_N) $ if $ u_j < v_j $ (i.e., $ u_j = 0 $, $ v_j = 1 $) at the first index $ j $ where $ u_j \neq v_j $. If multiple vertices yield the same minimum $ \chi_t^v / p_{\overrightarrow{H}}(u, v) $, select the smallest vertex $ v_* $ in lexicographic order. Set $ X_t^{(N)} = v_* $.
 
 It is immediate that conditional on the event $\{\ocube=\overrightarrow{H}\}$, the process $\mathbf{X}^{(N)}=(X^{(N)}_{t})_{t\in \Z_+}$ above constructed is a Markov chain with transition probability matrix $(p_{\overrightarrow{H}}(x,y))_{(x,y)\in \mathsf{vec}(\edges)}$.
Furthermore, the annealed law of $\mathbf{X}^{(N)}$ is given by \eqref{annealed} and the quenched law satisfies \cref{assump}.

\subsection{Absorbing sets clustering together}
For two vertices $ u = (u_1, u_2, \ldots, u_N)$ and $ v = (v_1, v_2, \ldots, v_N)$, we say that $u$  precedes $v$ in the \textit{lexicographic order} if $ u_j < v_j $ (i.e., $ u_j = 0 $, $ v_j = 1 $) at the first index $j$ where $ u_j\neq v_j$.

Let $(\eta^{v})_{v\in\vertices}$ be independent exponential random variables with parameter 1, which are also independent of {$(\chi_k^{\overrightarrow{e}})_{k\in \mathbb{Z}_+, \overrightarrow{e}\in \text{vec}(E_N)}$} and $(g_{i}(v))_{i\in [N],v\in \vertices}$. 

For each absorbing set $\mathcal{U} \in \Abso^*_N$, let $\des(\mathcal{U})$ denote a vertex in $\mathcal{U}$ defined as follows. If there exists a unique vertex $ v \in \mathcal{U} $ achieving the minimum value of $\eta^v$, set $$\des(\mathcal{U}) = \arg\min_{v \in \mathcal{U}} \eta^v.$$ If multiple vertices achieve the minimum, choose the smallest one in the lexicographic order. The latter case occurs with probability 0 as the $\eta^v$ are independent with continuous distributions. We note that by this construction, $\des(\mathcal{U})$ is uniformly chosen among vertices in $\mathcal{U}$.  We call $\des(\mathcal{U})$ the \textit{designated vertex} of $\mathcal{U}$. Let $\mathcal D$ be the set of designated vertices
$$
\mathcal D:=\{ u\in \vertices\colon\; u =\des(\mathcal{U}) \mbox{ for some }\mathcal{U}\in\Abso_N^*\}.
$$

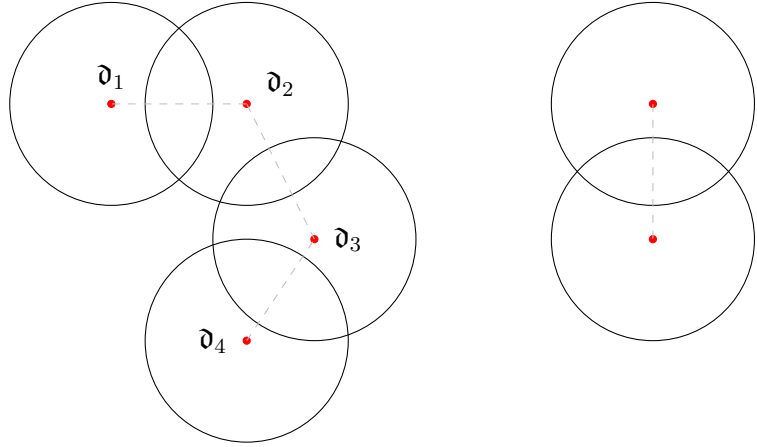
\begin{figure}[ht] \centering
\begin{tikzpicture}[scale=0.80]
\filldraw[color=black, fill=none](-2,1.5) circle (1.5);
\filldraw [red] (-2,1.5) circle (1.5pt);
\filldraw[color=black, fill=none](-1,3) circle (1.5);
\filldraw [red] (-1,3) circle (1.5pt);
\filldraw[color=black, fill=none](-2,5) circle (1.5);
\filldraw [red] (-2,5) circle (1.5pt);
\filldraw[color=black, fill=none](-4,5) circle (1.5);
\filldraw [red] (-4,5) circle (1.5pt);
\node at (-2.5,1.5) {$\mathfrak{d}_{4}$};
\node at (-0.5,3) {$\mathfrak{d}_{3}$};
\node at (-1.5,5.3) {$\mathfrak{d}_{2}$};
\node at (-4,5.4) {$\mathfrak{d}_{1}$};
\draw[lightgray] [dashed] (-2,1.5) -- (-1,3);
\draw [lightgray][dashed] (-1,3) -- (-2,5);
\draw [lightgray][dashed] (-2,5) -- (-4,5);
%\node at (-2,-0.5){};
\filldraw[color=black, fill=none](4,5) circle (1.5);
\filldraw [red] (4,5) circle (1.5pt);
\filldraw[color=black, fill=none](4,3) circle (1.5);
\filldraw [red] (4,3) circle (1.5pt);
\draw [lightgray][dashed] (4,5) -- (4,3);
%\node at (4,-0.5) {};
%\node at (4.5,5) {$v_{1}$};
%\node at (4.5,3) {$v_{2}$};
\end{tikzpicture}
\caption{\label{fi:cluster}{Every pair of designated vertices (in red) is connected with a dash line if the Hamming distance between them is less than or equal to $2m_{\beta}-2$. The union of discrete balls with radius $m_{\beta}-1$ centred at connected designated vertices forms a cluster.}}
\end{figure}

For each vertex $u \in \mathcal{V}_N$, let $\ll_u$ be an order on $\mathcal V_N$ which is defined as follows.
For two distinct vertices $x, y\in \mathcal{V}_N$, we denote $x \ll_u y$ if either $d(x,u)<d(y,u)$ or $d(x,u)=d(y,u)$ and $x$ is smaller than $y$ in the lexicographic order. Note that for any pair of distinct vertices $x$ and $y$, we have either $x \ll_u y$ or $y \ll_u x$. Set $$\Pi_N :=\{u\in \mathcal{V}_N: d(u,\mathcal{D})\le m_{\beta}-1\}.$$ 
Notice that $\Pi_N$ is not necessarily connected but it is partitioned into connected component(s). 
For each vertex $u\in \Pi_N$, let $V(u)$ be the connected component of $\Pi_N$ containing $u$. Let $\mathfrak{d}_1, \mathfrak{d}_2, \cdots, \mathfrak{d}_j$ be the designated vertices in $V(u)$ arranged in the ascending order corresponding to $\ll_u$. Set \begin{align}\label{def.Du}
    \mathcal{D}^{(u)}&:=\big\{\mathfrak{d}_k\ : \ 1\le k\le \min\{j, m_{\beta}\}\big\}\quad\text{and}\\ \label{def.Cu}
    \mathcal{C}(u)&:=\big\{v\in \mathcal{V}_N\ :\ d(v, \mathcal{D}^{(u)})\le m_{\beta}-1\big\}.\end{align}
We call $\mathcal{C}(u)$ the \textit{cluster} of $u$.  Each cluster covers a region of the hypercube, generated by at most $m_\beta $ designated vertices. %(hence at most $m_\beta$ absorbing sets, each of size at most $m_\beta$)
Roughly speaking, a cluster is a ``neighbourhood" of closely grouped absorbing sets. Notice that the process $\mathbf{X}^{(N)}$ visits clusters sequentially until it hits a trap or PNE. Within a cluster, local properties (e.g., presence of PNEs) can be studied independently. See \cref{fi:cluster} for a graphical example of two distinct clusters. 

We note that the process $\mathbf{X}^{(N)}$ may become trapped once it hits a designated vertex. For the sake of convenience, we extend the process to allow it to escape from a PNE or a trap as follows.

\begin{definition}\label{def.Xtilde}
Let $\sigma^{(N)}$ be the first time the process $\mathbf{X}^{(N)}$ hits a designated vertex, i.e.,
$$\sigma^{(N)} := \inf\{t \colon X^{(N)}_t \in \mathcal{D}\}.$$
We construct a new process $\widetilde{\mathbf{X}}^{(N)}=(\widetilde{X}_t^{(N)})_{t\in \mathbb{Z}_+}$ as follows. We set $\widetilde{X}_t^{(N)} := X_t^{(N)}$ for $0\le t\le \sigma^{(N)}$. Let $\widetilde{k}_{t}^{(u,v)}$ be the number of crosses by $\widetilde{\mathbf{X}}^{(N)}$ from $u$ to $v$ up to time $t$. For time $t>\sigma^{(N)}$ with $\widetilde{X}_{t-1}^{(N)}=u$, if there is more than one vertex $v\sim u$ achieving the same minimum $\chi^{(u,v)}_{\widetilde{k}_{t-1}^{(u,v)}}$, then we say the construction fails and the process remains at $u$ forever. Otherwise, we set
$$\widetilde{X}_{t}^{(N)}=\operatorname{argmin}_{v\sim u} \chi^{(u,v)}_{\widetilde{k}_{t-1}^{(u,v)}}.$$
Notice that the construction fails with probability 0. In other words, the process moves like a simple symmetric random walk on $\mathcal{H}_N$ after hitting $\mathcal{D}$. We call the process $\widetilde{\mathbf{X}}^{(N)}$ an \textit{extension} of $\mathbf{X}^{(N)}$.
\end{definition}

Note that $\sigma^{(N)}$ is a stopping time with respect to the filtration $(\mathcal{F}_n)_{n\in\Z_+}$ where $\mathcal{F}_n$ is generated by $(X_t)_{0\le t\le n}$ and $\mathcal{D}$. For the sake of simplicity, we use from now on the notations $\bf X$ and $\widetilde{\bf X}$ instead of ${\bf X}^{(N)}$ and $\widetilde{\bf X}^{(N)}$.

We define recursively a sequence of random times $(S_{k})_{k\in \Z_+}$  as follows. Set $S_0 = 0$, and let  
$$S_k := \inf\Big\{ t > S_{k-1}  \colon \widetilde{X}_t\in \Pi_N\setminus \cup_{j=1}^{k-1} \mathcal{C}{(\widetilde{X}_{S_j})}  \Big\}.$$
Note that $(S_k)_{k\in \N}$ are stopping times with respect to the filtration $(\widetilde{\mathcal{F}}_n)_{n\in\Z_+}$ where $\widetilde{\mathcal{F}}_n$ is generated by $(\widetilde{X}_t)_{1\le t\le n}$, $\ocube$ and $\mathcal{D}$.

We notice that $\mathcal{C}(\widetilde{X}_{S_k})\neq \mathcal{C}(\widetilde{X}_{S_j})$ for all $1\le j \le k-1$.
In other words, a new cluster $\mathcal{C}(\widetilde{X}_{S_k})$ is revealed by the process $\widetilde{\mathbf{X}}$ at each time $S_k$ for $k\in \N$.

 \begin{definition}\label{rm:mb2} We say that a cluster $\mathcal{C}(u)$ is \textit{good} if {$d(\mathcal{D}^{(u)}, \mathcal D\setminus \mathcal D^{(u)})\ge  2m_{\beta}$}.
We define the event $$\mathscr{J}_k :=\{S_k<\infty\}\cap \big\{ \mathcal{C}(\widetilde{X}_{S_k}) \text{ is good}\big\}.$$
     \end{definition}

In other words, a cluster is good if it is well separated from the other clusters. Recall from the definitions \eqref{def.Du}--\eqref{def.Cu} that the cluster $\mathcal{C}(u)$ is the union of balls of radius $m_\beta-1$ centred at the designated vertices in $\mathcal{D}^{(u)}$. {In the worst-case scenario, there is a trap of size $m_\beta$ in $\mathcal{D}^{(u)}$ and another trap of size $m_\beta$ in $\mathcal{D}\setminus\mathcal{D}^{(u)}$. A separation of at least $2m_\beta$ between $\mathcal{D}^{(u)}$ and $\mathcal{D}\setminus\mathcal{D}^{(u)}$ ensures that the orientations of the edges incident to the traps/PNEs containing the designated vertices in $\mathcal{D}^{(u)}$ are independent of those corresponding to $\mathcal{D}\setminus\mathcal{D}^{(u)}$.} Hence, good clusters make it possible to treat the behaviour of the process within each cluster as independent of that in the others.

Note that on $\mathscr{J}_k$, the cluster $\mathcal{C}(\widetilde{X}_{S_k})$ containing no more than $m_{\beta}$ designated vertices, and it is a connected component of $\Pi_N$.

 Let $$q_N:= \lfloor {\rm e}^{(\log N)^2}\rfloor.$$ {The following lemma asserts that, with high probability, all clusters $\mathcal{C}(\widetilde{X}_{S_k})$, for $1 \le k \le {q}_N$, are good.}

\begin{lemma} \label{lem.cluster}
    There exists a constant $\theta\in (0,1)$ such that for large enough $N$, 
    $$\Prob\Big(\bigcap_{k=1}^{q_N}\mathscr{J}_{k} \Big)\ge 1-\theta^{N}.$$
\end{lemma}  

\begin{proof} 
It is sufficient to prove that there exists a constant $\theta_0\in (0,1)$ such that for large enough $N$,
\begin{align} \label{prob.J}\Prob\Big(  {\Big(\bigcap_{i=1}^{k-1}\mathscr{J}_{i}\Big)}\cap \mathscr{J}_{k}^c \Big)\le \theta_0^{N}.
\end{align}
for $1\le k\le q_N$.

Fix an integer $k$ such that $1\le k\le q_N$. Recall that, by \cref{def.Xtilde}, the process $(\widetilde{X}_t)$ cannot be absorbed by any element in $\Abso_N^*$. On the event
$ \cap_{i=1}^{k-1} \mathscr{J}_{i}\cap \{S_k=\infty\}$ we have either the process is absorbed by an element $\mathcal{U}\in \Abso_N\setminus \Abso^*_N$ or $\mathcal {D} \subseteq \bigcup_{i=1}^{k-1} \mathcal{C}(\widetilde{X}_{S_i}).$ By reason of \cref{co:impoimpo},  the former case occurs with probability not larger than $\theta_1^N$ for some constant $\theta_1\in (0,1)$. In the latter case, the number of PNEs is not larger than $q_N m_{\beta}$ as each good cluster has at most $m_{\beta}$ designated vertices. Moreover, in virtue of \cite[Theorem 2]{amiet2019pure},
$$\mathbb P\big( \card{\pne_N} \le q_Nm_{\beta}\big)\le \theta_2^N,$$
for some $\theta_2\in (0,1)$. Hence
\begin{align}\label{case1}
   \Prob\left( {\Big(\bigcap_{i=1}^{k-1} \mathscr{J}_{i}\Big)}\cap \{S_k=\infty\}\right) \le \theta_1^N+\theta_2^N.
\end{align}
On the other hand, on the event $\{S_k<\infty\}\cap\big\{\mathcal{C}(\widetilde{X}_{S_k}) \text{ is not good}\big\}$, by \cref{rm:mb2}, there are more than $m_{\beta}$ designated vertices in $\ball_{4m^2_{\beta}}(\widetilde X_{S_k})$. Hence, by reason of \cref{le:McD2} and \cref{le: GinFrag}, we infer that 
\begin{align}\label{case2}
    \Prob(S_k<\infty,\ \mathcal{C}(\widetilde{X}_{S_k}) \text{ is not good}\big)\le \theta_3^N
\end{align}
for some constant $\theta_3\in (0,1)$.
Combining \eqref{case1} and \eqref{case2}, we obtain \eqref{prob.J}.
\end{proof}

In the remaining part of the paper, we use ${\rm Const}$ to denote a generic {\bf positive} constant that can vary in different expressions and is allowed to depend on $\alpha$, but does not depend on $N$. {For a subset $B\subset \vertices$, define
$$W(B)=\{x\in \vertices : d(x, B)\le m_{\beta}-1\}.$$
The next lemma states that, if a deterministic set of vertices is exactly the set of designated vertices in its $(m_\beta-1)$-neighbourhood, then, after conditioning on any information which does not reveal the local orientations deciding whether these vertices are PNEs or are designated vertices of traps, the set consists entirely of PNEs with high probability}

\begin{lemma}\label{lem.Bk} 
Let $\mathfrak d^*\subseteq\mathcal V_N$ be deterministic and satisfy
$1\le \card{\mathfrak d^*}\le m_\beta$. {Let $A$ be an event such that} \begin{itemize}
 
   \item   {$\Prob\left(A,\, W(\mathfrak d^*)\cap \mathcal D =\mathfrak d^*\right)>0.$,}
    \item {$A\in \mathcal{P}(\mathfrak{d}^*)$, where $\mathcal{P}(\mathfrak{d}^*)$ is the $\sigma$-algebra generated by
$(\chi_k^{\vec e})_{k\in\mathbb Z_+,\vec e\in
\mathsf{vec}(\edges)}$, $(\eta^v)_{v\in\mathcal V_N}$ and by the orientations of edges which are not incident
to any vertex belonging to an element of $\Abso_N^*$ whose designated vertex lies
in $W(\mathfrak d^*)$.} 
\end{itemize}
Then there exists a constant $\theta\in (0,1)$ such that, for all sufficiently large $N$,
\begin{equation*}
\Prob\left(
\mathfrak d^*\subseteq \pne_N
\mid {
A,\, W(\mathfrak d^*)\cap \mathcal D=\mathfrak d^*}
\right)
\ge 1-\theta^N.
\end{equation*}
\end{lemma}

\begin{proof}If $m_{\beta}<4$, then by the definition \eqref{Abso.star}, we have {$\Abso_N^*=\{\{v\}:v\in\pne_N\}$}. Therefore $\pne_N=\mathcal{D}$, and the lemma follows immediately. For the remainder of the proof, we assume that $m_{\beta}\ge 4$.

    Let $\mathfrak{d}_1, \mathfrak{d}_2, \ldots, \mathfrak{d}_j$ denote the distinct elements of $\mathfrak{d}^*$. Let $h$ denote the number of pairs of vertices in $\mathfrak{d}^*$ that are adjacent in the hypercube $\mathcal{H}_N$. Since $\card{\mathfrak{d}^*}\le m_{\beta}$, we have  $h \leq j^2 \leq m_{\beta}^2$. If the vertices $\mathfrak{d}_1, \mathfrak{d}_2, \ldots, \mathfrak{d}_j$ are all PNEs, then on the oriented hypercube $\ocube$, there are $h$ unoriented edges connecting the adjacent pairs in $\mathfrak{d}^*$, and $jN - 2h$ edges connecting $\mathfrak{d}^*$ and $\vertices\setminus \mathfrak{d}^*$, which are either unoriented or oriented from $\vertices \setminus \mathfrak{d}^*$ to $\mathfrak{d}^*$.
{Moreover, by the assumption on $A$, the event $A$ is independent of the orientations of the edges incident to $\mathfrak{d}^*$ if $\mathfrak{d}^*\subseteq\pne_N$.}
Hence
\begin{equation}\label{eq:bid}
 \Prob( \mathfrak{d}^*\subseteq\pne_N, {A})
= \alpha^{h}  \Big(\frac{1+\alpha}2\Big)^{j N-2h} {\Prob(A)}\ge \alpha^{m_{\beta}^2}  \Big(\frac{1+\alpha}2\Big)^{j N} {\Prob(A)}.
\end{equation}

We next show that
\begin{equation}\label{eq:final0}
\begin{aligned}
&\Prob\Big( W(\mathfrak d^*)\cap \mathcal D=\mathfrak d^*, \mathfrak{d}^*\not\subseteq \pne_N , {A} \Big)\\
&=\mathbb P\Big(W(\mathfrak d^*)\cap \mathcal D=\mathfrak d^*\; \mbox{ and $\mathfrak{d}^*\cap \mathcal{U}\neq \emptyset$ for some $\mathcal{U}\in {\sf Trap}^*_N$}, {A} \Big)\\
&\le  {\rm Const}\times  N^{m_\beta^2}\Big(\frac{1+\alpha}2\Big)^{(j+3) N}{\mathbb{P}(A).}
\end{aligned}
\end{equation}
The first identity follows immediately by the definition of designated vertices. Let {$R_1, R_2,\cdots, R_j$} be deterministic disjoint subsets of $\vertices$ which are possible outcomes of $j$ elements in $\Abso_N^*$ such that $\mathfrak{d}_k\in R_k$ for each $k\in\{1,2,\cdots,j\}$, and there exists $i\in \{1,2,\cdots, j\}$ such that $\card{R_i}\ge 4$, i.e. $R_i$ can be a vertex set of a trap. 
{This assumption implies that $1\le \card{R_k}\le m_\beta$ and $d(\mathfrak{d}_k, v)\le m_{\beta}-2$ for each $v\in R_k$, $1\le k\le j$.}
Let $\widetilde{h}$ be the number of pairs of vertices in $\bigcup_{k=1}^j R_k$ which are adjacent in the hypercube $\cube$.
If $R_1, R_2,\cdots, R_j$ are vertex sets in $\Abso_N^*$, then on $\ocube$, there are $N\sum_{k=1}^j \card{R_k} - 2\widetilde{h}$ edges connecting $\bigcup_{k=1}^j R_k$ and $\vertices\setminus \bigcup_{k=1}^j R_k$, which are either unoriented or oriented from  $\vertices\setminus\bigcup_{k=1}^j R_k$ to $\bigcup_{k=1}^j R_k$. 
{Moreover, by the assumption on $A$, $A$ is independent of the orientations of the edges incident to $\bigcup_{k=1}^j R_k$ if $R_1, R_2, \cdots, R_j$ are elements of $\Abso^*_N$.}
Therefore,
\begin{equation}\label{eq:event_subgp}
\begin{aligned}
	\Prob\left( W(\mathfrak d^*)\cap \mathcal D=\mathfrak d^*, \bigcap_{k=1}^j\{R_k \in \Abso^*_N\}, {A} \right) 
	&\le    \left(\frac{1+\alpha}{2}\right)^{ N\sum_{k=1}^j \card{R_k}-2\widetilde{h}}{\Prob(A)}\\
    & \le  \operatorname{Const}\times \Big(\frac{1+\alpha}2\Big)^{(j+3) N}{\Prob(A)}.
\end{aligned}
\end{equation}
where in the last inequality of \eqref{eq:event_subgp}, we used the following facts:
\begin{itemize}
    \item $\widetilde{h}\le m_{\beta}^4$ as there are at most $m_{\beta}^2$ vertices in $\bigcup_{k=1}^j R_k$. Here we recall that $j=\card{\mathfrak{d}^*}\le m_{\beta}$, and that $1\le \card{R_k}\le m_{\beta}$ for $1\le k\le j$ since $\mathbb  P(R_1,R_2,\cdots ,R_j\in \Abso^*_N)>0$.
    \item $\sum_{k =1}^j \card{R_k}\ge j+3$ as there exists $k\in \{1,2,\cdots, j\}$ such that $|R_k|\ge 4$.
\end{itemize}
For each $k\in \{1,2,\cdots, j\}$, we can define a tree rooted at $\mathfrak{d}_k$ with vertex set $R_k$. We note that, by virtue of Lemma 9 in \cite{McDScoWhi:2020}, the number of trees rooted at a vertex $v$ with $r$ vertices is bounded by $({\rm e} N)^{r-1}$. As $1\le \card{R_k}\le m_{\beta}$ and $j\le m_\beta$, we infer that the cardinality of all possible values of $R_1, R_2,\cdots, R_j$ is bounded by $\big(\sum_{r=1}^{m_{\beta}}({\rm e}N)^{r-1} \big)^{m_{\beta}}\le \text{Const}\times N^{m_{\beta}^2}$. Taking all possible values of $R_1, R_2,\cdots, R_j$ in \eqref{eq:event_subgp} and applying the law of total probability, we obtain \eqref{eq:final0}.

{We now prove that
\begin{equation}\label{eq:extraD}
\Prob\Big(\mathfrak d^*\subseteq \pne_N, W(\mathfrak d^*)\cap \mathcal D\neq \mathfrak d^*, A\Big)
\le
{\rm Const}\times N^{{2m_\beta}}
\Big(\frac{1+\alpha}{2}\Big)^{(j+1)N}\Prob(A).
\end{equation}
Indeed, on the event $\big\{\mathfrak d^*\subseteq \pne_N\big\}\cap\{ W(\mathfrak d^*)\cap \mathcal D \neq \mathfrak d^*\big\}\cap A$, there is an additional designated vertex $z\in W(\mathfrak d^*)\setminus \mathfrak d^*$. If $z$ is a PNE, then, by the same counting argument as in \eqref{eq:bid}, uniformly over $z\in W(\mathfrak d^*)$,}
$$
{
\Prob(\mathfrak d^*\cup\{z\}\subseteq\pne_N,A)
\le
{\rm Const}
\Big(\frac{1+\alpha}{2}\Big)^{(j+1)N}\Prob(A).
}
$$
{If $z$ belongs to an element $\mathcal U\in{\sf Trap}_N^*$, then the designated vertex of $\mathcal U$ is $z$ and $z\in W(\mathfrak d^*)$. Hence, by the assumption on $A$, the orientations of edges incident to $\mathcal U$ are independent of $A$. 
Since $\card{\mathcal U}\le m_\beta$, at least $N-m_\beta+1$ edges incident to $z$ connect $z$ to vertices outside $\mathcal U$, and all these edges are either unoriented or oriented towards $z$. Moreover, at most $m_\beta$ of these edges can also be involved in the event $\mathfrak d^*\subseteq\pne_N$, so the possible overlap changes the estimate only by a multiplicative constant. For fixed $z$, the number of possible connected sets $\mathcal U$ with $z\in\mathcal U$ and $\card{\mathcal U}\le m_\beta$ is bounded by $\sum_{r=1}^{m_\beta}({\rm e}N)^{r-1}\le {\rm Const}\times N^{m_\beta}$. Hence, uniformly over $z\in W(\mathfrak d^*)$,}
$$
{
\mathbb P\Big(
\mathfrak d^*\subseteq \pne_N,\ 
z\in \mathcal U \text{ for some } \mathcal U\in{\sf Trap}_N^*, A
\Big)
\le {\rm Const}\times N^{m_\beta}
\Big(\frac{1+\alpha}{2}\Big)^{(j+1)N}\Prob(A).
}
$$
{Since $\card{W(\mathfrak d^*)}\le {\rm Const}\times N^{m_\beta}$, \eqref{eq:extraD} follows by taking a union bound over all possible values of $z\in W(\mathfrak d^*)\setminus\mathfrak d^*$.}

{Combining \eqref{eq:bid} together with \eqref{eq:extraD} and using a union bound, we obtain that, for all sufficiently large $N$,}
\begin{equation}\label{eq:bid.exact}
{
\begin{aligned}
&\Prob\Big(
\mathfrak d^*\subseteq\pne_N,\,
W(\mathfrak d^*)\cap\mathcal D=\mathfrak d^*,\,
A
\Big)\\
&\ge
\Prob(\mathfrak d^*\subseteq\pne_N,A)
-
\Prob\Big(
\mathfrak d^*\subseteq\pne_N,\,
W(\mathfrak d^*)\cap\mathcal D\neq\mathfrak d^*,\,
A
\Big)\\
&\ge
\frac12
\alpha^{m_\beta^2}
\left(\frac{1+\alpha}{2}\right)^{jN}
\Prob(A).
\end{aligned}
}
\end{equation}

Notice that 
\begin{equation}
    {  \begin{aligned}\label{ineq.cond}
 & \Prob(\mathfrak{d}^*\subseteq \pne_N\mid  W(\mathfrak d^*)\cap \mathcal D=\mathfrak d^*, A)\\
    &
=\frac{\Prob(\mathfrak{d}^*\subseteq \pne_N,  W(\mathfrak d^*)\cap \mathcal D=\mathfrak d^*, A)}
{\Prob(\mathfrak{d}^*\subseteq \pne_N, W(\mathfrak d^*)\cap \mathcal D=\mathfrak d^*, A)
+\mathbb P\Big( W(\mathfrak d^*)\cap \mathcal D=\mathfrak d^*, \mathfrak{d}^*\not\subseteq \pne_N , A \Big)}.
\end{aligned}}
\end{equation}
Combining \eqref{ineq.cond} together with \eqref{eq:final0} and \eqref{eq:bid.exact}, we obtain
{\begin{align*}
  \Prob(\mathfrak{d}^*\subseteq \pne_N\mid W(\mathfrak d^*)\cap \mathcal D=\mathfrak d^* ,A)&\ge
\frac{
\frac12\alpha^{m_{\beta}^2}  \Big(\frac{1+\alpha}2\Big)^{j N}
}{
\frac12\alpha^{m_{\beta}^2}  \Big(\frac{1+\alpha}2\Big)^{j N}
+\operatorname{Const}\cdot N^{m_\beta^2}\Big(\frac{1+\alpha}2\Big)^{(j+3) N}
}\\
&\ge 1-\operatorname{Const} \cdot N^{m^2_\beta} \Big(\frac{1+\alpha}{2}\Big)^{3N},
\end{align*}}
yielding the result of the lemma.
\end{proof}

 We define \begin{align*}
 \mathscr{B}_k := \big\{\mathcal{D}\cap \mathcal{C}(\widetilde{X}_{S_k})  \subseteq \pne_N\big\}.\end{align*}
 
 Recall that  $q_N = \floor{{\rm e}^{(\log N)^2}}$. The following lemma asserts that, conditional on the event that all clusters $\mathcal{C}(\widetilde{X}_{S_k})$ for $1 \leq k \leq q_N$ are good, all the designated vertices in these clusters are PNEs with high probability.

\begin{lemma}\label{le:spe1}
There exists  $\theta \in (0,1)$ such that for all $N$ large enough we have that 
\begin{equation*}
 \Prob\left(\bigcap_{k=1}^{q_N} {\big(\mathscr{B}_k\cap \mathscr{J}_k\big)}\right)\ge 1- \theta^N.
\end{equation*}
\end{lemma}

\vspace{0.2cm}
\begin{proof}
Fix $k\in \{1,2,\cdots, q_N\}$. We first prove that there exists a constant $\theta_1\in (0,1)$ such that for $N$ large enough we have
\begin{align}\label{eq1}
    \Prob\Big(\mathscr{B}_k^c \cap  {\mathscr J}_k \Big) \le \theta_1^N.
\end{align}
Let $\mathfrak{d}^*\subseteq \vertices$ be a deterministic subset of vertices such that
\begin{equation}\label{eq:massth}
 \Prob\Big(\mathcal{C}(\widetilde{X}_{S_k}) \cap \mathcal{D} = \mathfrak{d}^*,{\mathscr J}_k\Big)>0.
\end{equation}
We recall that $$W(\mathfrak{d}^*):=\{u\in \mathcal{V}_N\ :\ d(u,\mathfrak{d}^*)\le m_{\beta}-1\}.$$
Note that \eqref{eq:massth} implies that $1\le \card{\mathfrak{d}^*}\le m_{\beta}$ and $W(\mathfrak{d}^*)$ is a possible outcome of a good cluster with the set of designated vertices $\mathfrak{d}^*$.
 Notice that on $\Big\{ \mathcal{C}(\widetilde{X}_{S_k}) \cap \mathcal{D} = \mathfrak{d}^*\Big\}\cap{\mathscr J}_k$ we have $\mathcal C(\widetilde X_{S_k})=W(\mathfrak d^*)$ and $W(\mathfrak d^*)\cap \mathcal D =\mathfrak d^*$. Set 
$$\mathscr{R}_k(\mathfrak d^*)
:=
\left\{
S_k<\infty,\
\mathcal C(\widetilde X_{S_k})=W(\mathfrak d^*),\
\mathcal C(\widetilde X_{S_k})\text{ is good}
\right\}.$$
Then
\begin{align}\label{rep.R}
\left\{
\mathcal C(\widetilde X_{S_k})\cap\mathcal D=\mathfrak d^*
\right\}
\cap
\mathscr J_k
= \mathscr R_k(\mathfrak d^*)
\cap
\left\{
 W(\mathfrak d^*)\cap \mathcal D=\mathfrak d^*
\right\}.
\end{align}

{Recall that $\mathcal{P}(\mathfrak d^*)$ is the $\sigma$-algebra generated by $(\eta^v)_{v\in \vertices}$, $(\chi_{k}^{\vec{e}})_{k\in \Z_+, \vec{e}\in \mathsf{vec}(\mathcal E_N)}$ and by the orientations of edges which are not incident to any vertex belonging to an element of $\Abso_N^*$ whose designated vertex lies in $W(\mathfrak d^*)$.}

{We next verify that, conditional on $\left\{
W(\mathfrak d^*)\cap\mathcal D=\mathfrak d^*
\right\}$, the event $\mathscr R_k(\mathfrak d^*)$ is determined by the information in $\mathcal{P}(\mathfrak d^*)$, i.e. 
\begin{align}\label{eq.sigma}
    \left\{
W(\mathfrak d^*)\cap\mathcal D=\mathfrak d^*
\right\} \cap  \mathscr R_k(\mathfrak d^*) \in \big\{ A\cap \{W(\mathfrak d^*)\cap\mathcal D=\mathfrak d^*\}  : A\in \mathcal{P}(\mathfrak d^*)  \big\}
\end{align}
Indeed, on $\left\{
W(\mathfrak d^*)\cap\mathcal D=\mathfrak d^*
\right\}$, the only designated vertices in $W(\mathfrak d^*)$ are the vertices of $\mathfrak d^*$.
For an element $\mathcal{U}\in \Abso_N^*$ with a designated vertex $\mathfrak{d}_i\in \mathfrak{d}^*$, we have $d(v, \mathfrak{d}_i)\le m_{\beta}-2$ for each $v\in \mathcal{U}$. Hence every edge incident to $\mathcal{U}$ has both endpoints in $W(\mathfrak d^*)$. On the event $\{W(\mathfrak d^*)\cap\mathcal D=\mathfrak d^*\} \cap  \mathscr R_k(\mathfrak d^*)$, $S_k$ is the first time at which the process enters the new cluster $W(\mathfrak d^*)$. Hence, the trajectory before $S_k$ does not use any edge whose two endpoints lie in $W(\mathfrak d^*)$. 
Therefore the orientations of edges incident to elements of $\Abso_N^*$ whose designated vertices lie in $W(\mathfrak d^*)$ are not needed to determine $\mathscr R_k(\mathfrak d^*)$. Hence \eqref{eq.sigma} is verified.}

{By \eqref{eq.sigma}, there exists an event $A_0\in \mathcal{P}(\mathfrak d^*)$ such that
\begin{align}\label{rep.A0}
    A_0
\cap
\left\{
W(\mathfrak d^*)\cap\mathcal D=\mathfrak d^*
\right\}
=
\mathscr R_k(\mathfrak d^*)
\cap
\left\{
W(\mathfrak d^*)\cap\mathcal D=\mathfrak d^*
\right\}.
\end{align}
Therefore, combining \eqref{rep.R}-\eqref{rep.A0} and applying \cref{lem.Bk} to the event $A_0$, we deduce that there exists $\theta_1\in (0,1)$ such that for all sufficiently large $N$,
\begin{align*} 
\begin{aligned}
&\Prob\left(
\mathfrak d^*\subseteq\pne_N
\mid
\mathcal C(\widetilde X_{S_k})\cap\mathcal D=\mathfrak d^*,
\mathscr J_k
\right)=
\Prob\left(
\mathfrak d^*\subseteq\pne_N
\mid
A_0,
 W(\mathfrak d^*)\cap \mathcal D =\mathfrak d^*
\right)\ge
1-\theta_1^N.
\end{aligned}
\end{align*}}
Since, on $\big\{
\mathcal C(\widetilde X_{S_k})\cap\mathcal D=\mathfrak d^*
\bigr\}$,
the event $\mathscr B_k$ is precisely $\{\mathfrak d^*\subseteq\pne_N\}$,
we obtain
$$\Prob\left(
\mathscr B_k^c
\mid
\mathcal C(\widetilde X_{S_k})\cap\mathcal D=\mathfrak d^*,
\mathscr J_k
\right)
\le \theta_1^N.$$
Using the law of total probability, where we take all possible values of $\mathfrak d^*$, we obtain \eqref{eq1}.

On the other hand, in view of \cref{lem.cluster}, there exists $\theta_2\in (0,1)$ such that for all sufficiently large $N$,
\begin{align}\label{ineq4}
 \Prob\Big(\bigcap_{k=1}^{q_N}  \mathscr{J}_k \Big)\ge 1- \theta_2^N. 
\end{align}

Finally, using a union bound, combining together with  \eqref{eq1} and \eqref{ineq4}, we get
\begin{align*}
 \Prob\left(\bigcap_{k=1}^{q_N} \big( \mathscr{B}_k\cap \mathscr{J}_k\big)\right)
    &\ge \Prob\Big(\bigcap_{k=1}^{q_N}  \mathscr{J}_k \Big) - \sum_{k=1}^{q_N}\Prob\big(\mathscr{B}_k^c \cap \mathscr{J}_k\big)\\
    &\ge 1-\theta_2^N-q_N\theta_1^N.
\end{align*}
This yields the result of the lemma.
\end{proof}

Let 
\begin{align*}
    \mathscr{A}_k=
\left\{
\widetilde X_{S_k+m_\beta-1}
\in
\mathcal C(\widetilde X_{S_k})\cap\mathcal D
\right\}
    \quad \text{and}\quad \mathscr{G}_k=\bigcap_{i=1}^{k}\big(\mathscr{A}_i^c\cap \mathscr{B}_i\cap \mathscr{J}_i\big).
\end{align*}

The event $\mathscr A_k$ means that, in exactly $m_\beta-1$ steps after entering the
$k$-th cluster, the process is at a designated vertex in that cluster. The
next lemma shows that, conditional on the events
\begin{itemize}[leftmargin=*, itemsep=0pt, topsep=0pt]
     \item $\bigcap_{i=1}^{k} \mathscr{J}_{i}$: for each $1 \le i \le k$,
    the process reaches the cluster $\mathcal C(\widetilde X_{S_i})$ in finite time, and that cluster is good;
    
    \item $\bigcap_{i=1}^{k} \mathscr{B}_i$: all designated vertices in these
    clusters are PNEs;
    
    \item $\bigcap_{i=1}^{k-1} \mathscr{A}_i^c$: for each $1\le i\le k-1$, in
    exactly $m_\beta-1$ steps after the process first enters the cluster
    $\mathcal C(\widetilde X_{S_i})$, the process is not at a designated vertex in that cluster;
\end{itemize}
the event $\mathscr{A}_k$ occurs with probability at least $K_1 N^{-K_2}$, for some positive constants $K_1, K_2$.

\begin{lemma}\label{le:spe2} There exist positive constants  $K_1$ and $K_2$ which {do not depend} on $N$ such that  for $1\le k\le q_N$, we have
$$\Prob( \mathscr{A}_k \;|\; \mathscr{B}_k \cap {\mathscr J}_k \cap \mathscr{G}_{k-1} ) \ge K_1 N^{- K_2}.$$
\end{lemma}

\begin{proof}
{The case $m_\beta=1$ is trivial. Hence for the remainder of the proof, we assume that $m_\beta\ge2$.} Fix $1\le k\le q_N$. Consider a deterministic subset $\mathfrak{d}^*\subseteq \mathcal{V}_N$ and a deterministic vertex $x_1\in \vertices$ such that  $\Prob(H)>0$, where
$$H:=\big\{\widetilde{X}_{S_k}=x_1\big\}\cap\big\{ \mathcal{C}(\widetilde{X}_{S_k})\cap \mathcal{D} =\mathfrak{d}^*\big\}\cap \mathscr{B}_k\cap \mathscr{J}_k \cap \mathscr{G}_{k-1}.$$
Note that on $H\subseteq \mathscr{J}_k \cap \big\{ \mathcal{C}(\widetilde{X}_{S_k})\cap \mathcal{D} =\mathfrak{d}^*\big\}$,  $\mathfrak{d}^*$ is the set of designated vertices in a good cluster. Since $\Prob(H)>0$, we must have $1\le \card{\mathfrak d^*}\le m_{\beta}$ and $d(x_1, \mathfrak{d}^*)= m_{\beta}-1$. 

Let  $x_1,x_2,\ldots,x_{m_{\beta}}$ be deterministic vertices on a shortest path in $\cube$ which connects $x_1$ and $\mathfrak{d}^*$ with $x_{m_{\beta}}\in \mathfrak{d}^*$. 
{Denote $e_j:= \{x_{j},x_{j+1}\}$ for $1\le j\le m_{\beta}-1$. To encode the orientation of $(e_j)_{1\le j\le m_{\beta}-1}$, we use random variables $Y:=(Y_j)_{1\le j\le m_{\beta}-1}$, where
$$Y_j:=\begin{cases}
1 & \text{if } (x_{j},x_{j+1})\in \oedges,\\
-1 & \text{if } (x_{j+1},x_{j})\in \oedges,\\
0 & \text{otherwise}.
\end{cases}$$}

{The random variable $Y$ takes values in $S:=\{-1,0,1\}^{m_{\beta}-1}$. We now apply \cref{lem:mono} (see in the appendix) to the random variable $Y$ and the event $H$.
Let $\mathcal G$ be the $\sigma$-algebra generated by $(\eta^v)_{v\in \vertices}$, $(\chi_{k}^{\vec{e}})_{k\in \Z_+, \vec{e}\in \mathsf{vec}(\mathcal E_N)}$ and by the orientations of all edges except $e_1,\cdots, e_{m_{\beta}-1}$. Thus $Y$ is independent of $\mathcal{G}$ and $H\in \mathcal G\vee\sigma(Y)$. 
For each $s\in S$, define
 $$H_s:=\big\{\mathbb{E}[\mathbf{1}_H\mathbf{1}_{\{Y=s\}}\mid\mathcal{G}] =\Prob(Y=s)\big\}\in \mathcal{G}.$$
By \cref{lem:mono}, for each $s\in S$, we have a.s. 
\begin{align}\label{eq.proj}
H_s\cap\{Y=s\}=
H\cap\{Y=s\}.
\end{align} 
Roughly speaking, $H_s$ is the event, depending only on the randomness in $\mathcal G$, on which the event $H$ occurs when the orientations of $e_1,\ldots,e_{m_\beta-1}$ encoded by $Y$ are fixed to be $s$.}

{
Set $s_*:=(1,1,\cdots, 1)\in S$. We next verify that a.s.
\begin{align}\label{H.incl}
    H_s\subseteq H_{s_*} \quad \text{for all }s\in S.
\end{align}
Fix $s\in S$. Since $H_s\in\mathcal G$, the occurrence of $H_s$ is determined by the variables in $\mathcal G$, and hence is unaffected by the orientations of $e_1,\ldots,e_{m_\beta-1}$. Moreover, by \eqref{eq.proj}, whenever a realisation of the variables in $\mathcal G$ belongs to $H_s$, the event $H$ occurs if the orientations of $e_1,\ldots,e_{m_\beta-1}$ are chosen so that $Y=s$. Thus, to prove \eqref{H.incl}, it is enough to show that for such a realisation of the variables in $\mathcal G$, the event $H$ still occurs when all variables in $\mathcal G$ are kept fixed and only the orientations of $e_1,\ldots,e_{m_\beta-1}$ are chosen so that $Y=s_*$. In other words, it is enough to verify that, under this replacement, the events
$\{\widetilde X_{S_k}=x_1\}$, $\{\mathcal C(\widetilde X_{S_k})\cap\mathcal D=\mathfrak d^*\}$, $\mathscr B_k$, $\mathscr J_k$, and $\mathscr G_{k-1}$
still hold. This claim follows from the following observations:
\begin{itemize}[leftmargin=*, itemsep=0pt, topsep=0pt]
    \item Replacing $Y$ by $s_*$ only changes the orientations of the edges $e_j$, with $1\le j\le m_\beta-1$, and does not change the trajectory of $\widetilde{\mathbf{X}}$ before time $S_k$. Indeed, all the edges $e_j$, with $1\le j\le m_{\beta}-1$, have both endpoints in $\mathcal C(\widetilde X_{S_k})$, and before time $S_k$, the process has not entered this cluster. Hence $\{\widetilde X_{S_k}=x_1\}$ still holds under this replacement.
    \item We note that every vertex in $\mathfrak d^*$ remains a PNE after changing $Y$ to $s_*$. Indeed, on $H_s\cap\{Y=s\}$, we have $\mathfrak d^*\subseteq\pne_N$. Changing $Y$ only changes the orientations of the edges $e_1,\ldots,e_{m_\beta-1}$. Since no vertex among $x_1,\ldots,x_{m_\beta-1}$ belongs to $\mathfrak d^*$, the only vertex of $\mathfrak d^*$ which may be incident to a changed edge is $x_{m_\beta}$. Under the replacement $Y=s_*$, the edge $\{x_{m_\beta-1},x_{m_\beta}\}$ is oriented toward $x_{m_\beta}$. Hence this change is compatible with $x_{m_\beta}$ remaining a PNE, and all other vertices of $\mathfrak d^*$ are unaffected. Therefore $\mathfrak d^*\subseteq\pne_N$ still holds when we change $Y$ to $s_*$. 
    \item Under the replacement $Y=s_*$, no vertex among $x_1,\dots,x_{m_{\beta}-1}$ belongs to $\mathcal{D}$. Indeed, on this event, for $1\le \ell \le m_{\beta}-1$, we have $(x_\ell,x_{\ell+1})\in\oedges$ and therefore, $x_{\ell}$ cannot be a PNE. If $x_{\ell}$  is a vertex of a trap for some $1\le \ell \le m_{\beta}-1$, then $x_{m_{\beta}}$ must be also a vertex of the same trap, which contradicts the fact from the previous item that $x_{m_\beta}$ remains a PNE after changing $Y$ to $s_*$.  Therefore no new designated vertex is created in $\mathcal C(\widetilde X_{S_k})$ when changing $Y=s_*$.
    Hence $\{\mathcal C(\widetilde X_{S_k})\cap\mathcal D=\mathfrak d^*\}$ still holds under this replacement.
    \item Note that, on the event $\{\mathcal C(\widetilde X_{S_k})\cap\mathcal D=\mathfrak d^*\}$, the event $\mathscr B_k$ is equivalent to $\{\mathfrak d^*\subseteq\pne_N\}$. 
    Hence $\mathscr{B}_k$ still holds when we change $Y$ to $s_*$.
    \item Since changing $Y=s_*$ does not create new designated vertices,  the cluster $\mathcal{C}(\widetilde{X}_{S_k})$ remains good. Therefore, the event $\mathscr J_k$ still holds under this replacement..
\item For $1\le i\le k-1$, replacing $Y$ by $s_*$ only changes the orientations of edges contained in the $k$-th cluster, and hence the events $\mathscr B_i$ and $\mathscr J_i$ remain unchanged. Since the $i$-th and $k$-th clusters are good and distinct, we have $d(\mathcal C(\widetilde X_{S_i})\cap\mathcal D,\mathfrak d^*)\ge 2m_\beta$, while $d(x_1,\mathfrak d^*)=m_\beta-1$; hence $d(x_1,\mathcal C(\widetilde X_{S_i})\cap\mathcal D)\ge m_\beta+1$. If $S_i+m_\beta-1<S_k$, then $\mathscr A_i^c$ is determined before time $S_k$, and so it is unchanged. If $S_i+m_\beta-1\ge S_k$, then $S_i+m_\beta-1-S_k\le m_\beta-2$, and starting from $x_1$ at time $S_k$, the process cannot reach $\mathcal C(\widetilde X_{S_i})\cap\mathcal D$ by time $S_i+m_\beta-1$. Therefore $\mathscr A_i^c$ remains true for every $1\le i\le k-1$, and consequently $\mathscr G_{k-1}=\bigcap_{i=1}^{k-1}(\mathscr A_i^c\cap \mathscr B_i\cap\mathscr J_i)$ still holds when we change $Y$ to $s_*$.
\end{itemize}
Hence \eqref{H.incl} is verified. Applying \cref{lem:mono} to the event $H$ and the random variable $Y$, we deduce that
\begin{align*}\label{eq:H.mono}
\Prob(Y=s_*\mid H)\ge \Prob(Y=s_*).
\end{align*}
Since $\Prob(Y=s_*)=\beta^{m_{\beta}-1}$, it follows that}
\begin{equation}\label{eq.path1}
    \begin{aligned}
       &\Prob\big( (x_{j},x_{j+1})\in \oedges \text{ for all } 1\le j \le m_{\beta}-1 \mid H \big)\ge \beta^{m_{\beta}-1}. 
    \end{aligned}
\end{equation}

Using \eqref{eq.path1} and \cref{assump}, we infer that for all sufficiently large $N$, 
 \begin{align*}
 & \Prob\Big( \mathscr{A}_{k} \;|\; \widetilde{X}_{S_k}=x_1, \mathcal{C}(\widetilde{X}_{S_k})\cap \mathcal{D} =\mathfrak{d}^*, \mathscr{B}_k, \mathscr{J}_k , \mathscr{G}_{k-1} \Big)\\
 & \ge \Prob\Big(\widetilde{X}_{S_k+j}=x_{j+1} \text{ for each } 1\le j \le m_{\beta}-1 \;\big|\; \widetilde{X}_{S_k}=x_1, \mathcal{C}(\widetilde{X}_{S_k})\cap \mathcal{D} =\mathfrak{d}^*, \mathscr{B}_k, \mathscr{J}_k , \mathscr{G}_{k-1} \Big) \\
 &\ge  \beta^{m_{\beta}-1}(\kappa_1N^{-\kappa_2})^{m_{\beta}-1}, 
 \end{align*}
 in which the last inequality follows from the strong Markov property and \cref{assump}. The proof ends by using the law of total probability over $\mathfrak{d}^*$ and $x_1$.
\end{proof}

 We now proceed to prove our main result.
 \begin{proof}[Proof of \cref{th:BRD}] Recall that $q_N = \floor{{\rm e}^{(\log N)^2}}$.
In view of \cref{le:spe1}, there exists $\theta\in (0,1)$ such that, for sufficiently large $N$,
\begin{equation}
\Prob\left(\bigcap_{k=1}^{q_N} \big( \mathscr{B}_k\cap {\mathscr J}_{k}\big) \right)\ge {1-  \theta^N.}
\end{equation}
Using \cref{le:spe2}, we obtain that, for sufficiently large $N$,
\begin{align}
\nonumber \Prob\left(\bigcap_{k=1}^{q_N}\big(\mathscr{A}_k^c\cap  \mathscr{B}_k\cap \mathscr{J}_{k}\big) \right)
 &\le \left(1- K_1 N^{- K_2}\right)^{q_N}\\
\label{eq:righsu} &\le {\rm e}^{- K_1 q_N N^{- K_2}}.
 \end{align}
The right-hand side of \eqref{eq:righsu} is summable. Hence $\big(\bigcup_{k=1}^{q_N} \mathscr{A}_k\big) \cap \big( \bigcap_{k=1}^{q_N}  \mathscr{B}_k\cap {\mathscr J}_{k}\big)$ occurs with high probability.
We notice that
\begin{align*}
\Big(\bigcup_{k=1}^{q_N} \mathscr{A}_k\Big) \cap \Big(\bigcap_{k=1}^{q_N} \big( \mathscr{B}_k\cap {\mathscr J}_{k}\big)\Big)
&\subseteq \bigcup_{k=1}^{q_N}\Big( \mathscr{A}_k \cap \bigcap_{i=1}^{k-1} \mathscr{A}_{i}^c \cap \bigcap_{i=1}^{k} \big(\mathscr{B}_i \cap {\mathscr J}_{i}\big)\Big)\\
&= \bigcup_{k=1}^{q_N} \big(\mathscr{A}_k\cap \mathscr{B}_k\cap \mathscr{J}_k \cap \mathscr{G}_{k-1}\big),
\end{align*}
where we recall that  $\mathscr{G}_{k}=\bigcap_{i=1}^{k}\big(\mathscr{A}_i^c\cap \mathscr{B}_i\cap \mathscr{J}_i\big)$. It follows that with high probability there exists $k\in [q_N]$ such that $\mathscr{A}_k\cap \mathscr{B}_k\cap \mathscr{J}_k \cap \mathscr{G}_{k-1}$ occurs.

Recall that $\sigma^{(N)} := \inf\{t\ge 0 \colon \widetilde{X}^{(N)}_t \in \Des\}$ and that $\widetilde X_t=X_t$ for all $0\le t\le \sigma^{(N)}$. 
{We now work on the event $\mathscr{A}_{k}\cap \mathscr{B}_{k}\cap\mathscr{J}_{k}\cap  \mathscr{G}_{k-1}$. On this event, by the definition of $\mathscr A_k$, we have $\sigma^{(N)}\le S_k+m_\beta-1$. We claim that  $\widetilde X_{\sigma^{(N)}}$ is a PNE on $\mathscr{A}_{k}\cap \mathscr{B}_{k}\cap\mathscr{J}_{k}\cap  \mathscr{G}_{k-1}$. Indeed, since every designated vertex belongs to $\Pi_N$, and since the $k$-th cluster is good, the first time the process enters the cluster containing $\widetilde X_{\sigma^{(N)}}$ is one of the times $S_i$ with $1\le i\le k$. Therefore, by the definition of $\mathscr B_i$ for the corresponding cluster, the designated vertex $\widetilde X_{\sigma^{(N)}}$ is a PNE. Consequently, 
$$\mathscr{A}_{k}\cap \mathscr{B}_{k}\cap\mathscr{J}_{k}\cap  \mathscr{G}_{k-1}\subseteq \big\{\tau^{(N)}=\sigma^{(N)}\big\}.$$
We next show that, on $\mathscr{A}_{k}\cap \mathscr{B}_{k}\cap\mathscr{J}_{k}\cap  \mathscr{G}_{k-1} \cap\{\Abso_N=\Abso_N^*\}$, no trap is hit before time $\sigma^{(N)}$. Indeed, on the event $\{\Abso_N=\Abso_N^*\}$, every trap has at most $m_\beta$ vertices, and thus every vertex of a trap is within distance at most $m_\beta-1$ of its designated vertex. Hence every trap vertex belongs to $\Pi_N$ and lies in the cluster containing the designated vertex of this trap. Thus, if the process hit a trap before time $\sigma^{(N)}$, then it would have entered the cluster containing the designated vertex of this trap before time $\sigma^{(N)}\le S_k+m_\beta-1$. The corresponding discovery time is some $S_i$ with $1\le i\le k$. However, by the definition of $\mathscr B_i$, all designated vertices in this cluster are PNEs, which is incompatible with one of them being a designated vertex of a trap. Hence no trap is hit before time $\sigma^{(N)}$. Consequently,}
$$
\mathscr{A}_{k}\cap \mathscr{B}_{k}\cap\mathscr{J}_{k}\cap  \mathscr{G}_{k-1} \cap\{\Abso_N=\Abso_N^*\}
\subseteq
\{ \tau^{(N)}<\xi^{(N)}\}.
$$

In view of \cref{co:impoimpo}, we note that $\Abso_N=\Abso_N^*$ with high probability. Hence $\tau^{(N)}<\xi^{(N)}$ with high probability.
\end{proof}

    \section*{Appendix}

 \begin{lemma}\label{lem:mono}{
Let $(\Omega,\mathcal F,\Prob)$ be a probability space. Let $\mathcal G$ be a
sub-$\sigma$-algebra of $\mathcal F$, and let $Y$ be a random variable taking
values in a finite set $S$. Assume that $Y$ is independent of $\mathcal G$ and $\Prob(Y=s)>0$ for every $s\in S$. Let $H\in \mathcal G\vee\sigma(Y)$ be an event such that $\Prob(H)>0$. For each $s\in S$, define}
$${
H_s :=
\left\{
\mathbb E\left[\mathbf 1_H\mathbf 1_{\{Y=s\}}\mid\mathcal G\right]
=
\Prob(Y=s)
\right\}.}
$$
\begin{itemize}
    \item[(i)] {We have a.s. $H_s\cap \{Y=s\}=H\cap \{Y=s\}$ for each $s\in S$.}
    \item[(ii)] {Assume that there exists $s_*\in S$ such that a.s. $H_s\subseteq H_{s_*}$ for each $s\in S$.
Then}
\begin{align*}{
\Prob(Y=s_*\mid H)\ge \Prob(Y=s_*).}
\end{align*}
\end{itemize}
\end{lemma}

\begin{proof}{
By the definition of $H_s$, we have $H_s\in\mathcal G$ for each $s\in S$. Since $H\in \mathcal G\vee\sigma(Y)$ and $S$ is finite, for each $s\in S$ there exists  $A_s\in \mathcal G$ such that 
$H\cap\{Y=s\}=A_s\cap\{Y=s\}$. We show that a.s. $A_s=H_s$. Indeed, since $Y$ is independent of $\mathcal{G}$, we have a.s.
$$\mathbb E[\mathbf 1_H\mathbf 1_{\{Y=s\}}\mid \mathcal G]
= \mathbb E[\mathbf 1_{A_s}\mathbf 1_{\{Y=s\}}\mid \mathcal G] =
\mathbf 1_{A_s}\Prob(Y=s).$$ 
Since $\Prob(Y=s)>0$, we thus have
$
H_s
=
\{\mathbf 1_{A_s}\Prob(Y=s)=\Prob(Y=s)\}
=
\{\mathbf 1_{A_s}=1\}
=
A_s$ a.s. Hence (i) is verified.}

We now prove (ii). Since $S$ is finite, it follows from (i) that a.s.
$
   H=\bigcup_{s\in S}\big(H_s\cap\{Y=s\}\big), 
$
where the union is disjoint. Therefore
\begin{align*}
\Prob(H)
=
\sum_{s\in S}\Prob\bigl(H_s\cap\{Y=s\}\bigr)=\sum_{s\in S}\Prob(H_s)\Prob(Y=s).
\end{align*}
where in the last step, we use the fact that $H_s\in\mathcal G$ for each $s\in S$ and $Y$ is independent of
$\mathcal G$.
Since a.s. $H_s\subseteq H_{s_*}$ for every $s\in S$, we have
$
\Prob(H_s)\le \Prob(H_{s_*})$
for every $s\in S.$
Therefore
\begin{align*}
\Prob(H)
&=
\sum_{s\in S}\Prob(H_s)\Prob(Y=s) \le
\sum_{s\in S}\Prob(H_{s_*})\Prob(Y=s) =
\Prob(H_{s_*}).
\end{align*}
Again using (i) and the independence of $Y$ and $\mathcal G$, we obtain
\begin{align*}
\Prob\bigl(H\cap\{Y=s_*\}\bigr)
&=
\Prob\bigl(H_{s_*}\cap\{Y=s_*\}\bigr)=
\Prob(H_{s_*})\Prob(Y=s_*).
\end{align*}
Since $\Prob(H)>0$, it follows that
\begin{align*}
\Prob(Y=s_*\mid H)
&=
\frac{\Prob\bigl(H\cap\{Y=s_*\}\bigr)}{\Prob(H)} 
=
\frac{\Prob(H_{s_*})\Prob(Y=s_*)}{\Prob(H)}
\ge
\Prob(Y=s_*),
\end{align*}
where the last inequality follows from $\Prob(H)\le \Prob(H_{s_*})$.
This completes the proof.
\end{proof}

\section*{Acknowledgement} We would like to thank   Alex Scott and  Marco Scarsini for very helpful discussion. The work of AC and TN was partially supported by ARC DP230102209. ZZ was supported by the Australian Government Research Training Program (RTP) Scholarship.

\bibliography{bibNEpercolation.bib}
\bibliographystyle{siam}

\end{document}